\documentclass[10pt]{amsart}

\usepackage{amsfonts}
\usepackage{amssymb}
\theoremstyle{plain}
\newtheorem{thm}{Theorem}[section]
\newtheorem{prop}[thm]{Proposition}
\newtheorem{cor}[thm]{Corollary}
\newtheorem{lemma}[thm]{Lemma}

\renewcommand{\mod}[1]{\left\lvert #1 \right\rvert}
\newcommand{\inp}[2]{\left\langle #1,#2 \right\rangle}
\newcommand{\norm}[1]{\left\| #1 \right\|}
\newcommand{\cl}[1]{\overline{#1}}
\newcommand{\bb}[1]{\mathbb{#1}}
\newcommand{\cal}[1]{\mathcal{#1}}
\newcommand{\cycgroup}[1]{\left\langle #1\right\rangle}

\DeclareMathOperator{\lat}{Lat}

\DeclareMathOperator{\alg}{Alg}

\DeclareMathOperator{\linspan}{span}
\DeclareMathOperator{\fix}{Fix}

\newcommand{\wk}{\text{weak$^\ast$}}
\DeclareMathOperator*{\mult}{mult}

\allowdisplaybreaks

\begin{document}

\title[Interpolation in $H^\infty_\Gamma$]{Abrahamse's interpolation theorem \\ and Fuchsian groups}

\author[Mrinal Raghupathi]{Mrinal Raghupathi}
\date{\today}
\address{Department of Mathematics, University of Houston \\ Houston, Texas 77204-3476, U.S.A.}
\email{mrinal@math.uh.edu}
\urladdr{http://www.math.uh.edu/~mrinal}
\thanks{The author wishes to acknowledge the support and guidance of his advisor, Vern I. Paulsen. This work was completed as part of the author's Ph.D. thesis at the University of Houston in August 2008.}
\begin{abstract}
We generalize Abrahamse's interpolation theorem from the setting of a multiply connected domain to that of a more general Riemann surface. Our main result provides the scalar-valued interpolation theorem for the fixed-point subalgebra of $H^\infty$ associated to the action of a Fuchsian group. We rely on two results from a paper of Forelli. This allows us to prove the interpolation result using duality techniques that parallel Sarason's approach to the interpolation problem for $H^\infty$. In this process we prove a more general distance formula, very much like Nehari's theorem, and obtain relations between the kernel function for the character automorphic Hardy spaces and the Szeg\"o kernel for the disk. Finally, we examine our interpolation results in the context of the two simplest examples of Fuchsian groups acting on the disk.
\end{abstract}

\maketitle

%\tableofcontents

\section{Introduction}\label{introduction}

\subsection{Motivation}
Our objective in this paper is a generalization of Abrahamse's interpolation theorem~\cite[Theorem~1]{abrahamse1}. Let $R\subseteq \bb{C}$ be a bounded, multiply connected domain of genus $g$. Let $H^\infty(R)$ be the Banach algebra of bounded holomorphic functions on $R$. Abrahamse considered a natural class of Hilbert modules $\{H^2_\lambda(R)\,:\,\lambda\in \bb{T}^g\}$ over the algebra $H^\infty(R)$. Each of the spaces $H^2_\lambda(R)$ is a reproducing kernel Hilbert space on $R$ with kernel function $K^\lambda$. Abrahamse's solution to the Pick interpolation theorem for the region $R$ is the following: 
\begin{thm}[Abrahamse]
Let $z_1,\ldots,z_n\in R$ and $w_1,\ldots,w_n\in \bb{D}$. Then there exists a holomorphic function $f:R\to \bb{D}$ such that $f(z_j) = w_j$ if and only if the matrices 
\[A_\lambda :=\left[(1-w_i\cl{w_j})K^\lambda(z_i,z_j)\right]\geq 0,\]
for all $\lambda\in \bb{T}^g$.
\end{thm}

There has been a considerable amount of  work related to Abrahamse's theorem: see Ball~\cite{ball1}, for the matrix-valued generalization; Fedorov and Vinnikov~\cite{fedorovvinnikov}, for a careful analysis of the scalar-valued problem; and McCullough~\cite{mccullough}, and McCullough and Paulsen~\cite{mcculloughpaulsen}, for the operator algebraic viewpoint. 

In this paper we take a viewpoint that in some sense predates all the above proofs. Let $p:\bb{D}\to R$ be the universal covering map. Associated to this covering map is the \textit{group of deck transformations}, that is, the group $\Gamma$ of M\"obius transformations such that $p\circ \gamma = p$ for all $\gamma\in\Gamma$. Note that the map $p^\ast:H^\infty(R)\to H^\infty(\bb{D})$ given by $p^\ast (f) = f\circ p$ is an isometric embedding of $H^\infty(R)$, and the range of $p^\ast$ is precisely the fixed-point subalgebra for the action of $\Gamma$, that is, the set of functions $f\in H^\infty(\bb{D})$ such that $f\circ\gamma = f$ for all $\gamma\in \Gamma$. Therefore, function theoretic problems on the domain $R$ can be lifted to the disk. This idea is a fundamental tool in the theory of Riemann surfaces~\cite{farkaskra}. If $\Gamma$ is the group of deck transformations associated to $p$, then $\Gamma$ is a discrete subgroup of the group of all M\"obius transformations on the disk. In addition, $\Gamma$ is torsion-free, $\Gamma$ acts without fixed points and the action is properly discontinuous. Further, the quotient space $\bb{D} \backslash \Gamma$ is a Riemann surface that can be identified with $R$.

A discrete group of M\"obius transformations is called a \textit{Fuchsian group}. Let $\Gamma$ be a Fuchsian group. We will consider the analogue of Abrahamse's theorem for $H^\infty_\Gamma$, the fixed-point subalgebra of $H^\infty$. 

Our approach in this paper owes a great deal to the elegant results obtained by Forelli~\cite{forelli}. Forelli initiated the study of function theory in these fixed-point algebras, and in the special case where $\Gamma$ is the group of deck transformations, he constructs a bounded projection from $P:H^\infty\to H^\infty_\Gamma$ that is also a bimodule map over $H^\infty_\Gamma$. The projection is obtained by modifying the natural conditional expectation from $H^\infty$ onto $H^\infty_\Gamma$. An explicit formula for the conditional expectation was obtained by Earle and Marden~\cite{earlemarden} and this was generalized by Ball~\cite{ball}. The construction of the expectation in these papers makes use of Poincar\'e series. The approach to function theory on multiply connected domains through the use of holomorphic vector bundles was initiated by Abrahamse and Douglas~\cite{abrahamsedouglas}. Their paper makes use of the Forelli projection as well. We will need neither the full strength of Forelli's results nor the extensions obtained in~\cite{ball} and~\cite{earlemarden}. Our proof rests on two lemmas in~\cite{forelli}. 

\subsection{Notation}
We will denote the open unit disk in the complex plane by $\bb{D}$. The circle will be denoted $\bb{T}$ and the $m$-torus will be denoted $\bb{T}^m$. By an \textit{automorphism} we mean a M\"obius transformation, that is, a holomorphic map of $\bb{D}$ onto $\bb{D}$ with a holomorphic inverse. An automorphism is the composition of a simple Blaschke factor and a rotation, that is, an automorphism $\phi$ is of the form 
\[ \phi(z)= \lambda \frac{a-z}{1-\cl{a}z},\]
where $a\in \bb{D}$ and $\lambda\in \bb{T}$. 
The derivative $\phi$ is easily computed to be 
\[\phi'(z) = \lambda \frac{\mod{a}^2-1}{(1-\cl{a}z)^2}.\]
Hence, given an automorphism $\phi$, the constants $a$ and $\lambda$ are determined by the equations
\[a = \phi^{-1}(0), \lambda = (\mod{a}^2-1)^{-1}\phi'(0).\]
The set of all automorphisms of the disk is a group under the composition of maps. This group is naturally identified with $PSL(2,\mathbb{R})$.

We will use the term \textit{automorphism group} to mean a subgroup of the group of all automorphisms of the disk and we will usually denote an automorphism group by $\Gamma$. An automorphism of the disk is a well-defined holomorphic map on the open set $\{z\,:\,\mod{z}<\mod{a}^{-1}\}\supseteq \bb{D}\cup\bb{T}$. Therefore, the action of $\Gamma$ on the disk extends to the boundary circle and defines an action of $\Gamma$ on the circle $\bb{T}$. 

Given a group $\Gamma$ of automorphisms and a point $w\in \bb{D}$ we denote the stabilizer of the point $w$ by $\fix(\Gamma,w)$, that is, 
\[\fix(\Gamma,w) := \{\gamma\in \Gamma\,:\, \gamma(w)=w\}.\] 
We denote the orbit of the point $w$ under the action of $\Gamma$ by $\Gamma(w)$, that is, 
\[\Gamma(w) := \{\gamma(w)\,:\, \gamma\in \Gamma\}.\] 
We will call two groups $\Gamma_1$ and $\Gamma_2$ \textit{conjugate} if and only if there exists an automorphism $\phi$ such that $\phi\Gamma_1 = \Gamma_2\phi$. We denote by $[\Gamma,\Gamma]$ the commutator subgroup of $\Gamma$.

The group $\Gamma$ acts, by composition, on the Lebesgue spaces $L^p:=L^p(\bb{T})$, for $1\leq p\leq\infty$, with respect to normalized Lebesgue measure $m$. Since a M\"obius transformation is analytic this action restricts to the Hardy spaces $H^p$.

If $X\subseteq L^p$ we denote by $[X]_p$ the smallest closed subspace of $L^p$ that contains $X$. If $p=2$, then we write $[X]_2 = [X]$. When $p=\infty$ we take closures with respect to the \wk{} topology.

The character group, or dual group, of $\Gamma$ is the set of continuous homomorphisms from the group $\Gamma$ into the circle $\bb{T}$. Let $X\subseteq L^p$ be a set that is closed under composition by elements of $\Gamma$. If $\sigma\in \hat{\Gamma}$, then we define 
\[X_\sigma:=\{f\in X\,:\, f\circ \gamma = \sigma(\gamma) f\}.\] 
We call this the \textit{character automorphic space associated to $\sigma$}. We call elements of $X_\sigma$ \textit{character automorphic}.

\subsection{Outline of our work}
It is fair to say that one of the best understood action is the irrational rotation. However, the elements of $L^1$ fixed by this action are constant. Since our focus is on function theoretic results, we are interested in the cases where the group fixes some nonconstant function in $H^\infty_\Gamma$, that is, when $H^\infty_\Gamma$ is nontrivial. If $f\in H^\infty_\Gamma$ is non trivial, then we can, after subtracting a constant, assume that $f(0)=0$. Therefore, $f(\gamma(0))=0$ for all $\gamma\in \Gamma$ and so the points in the orbit of the origin $\Gamma(0)$ must satisfy the Blaschke condition
\begin{equation}\label{blaschkesum}
\sum_{\zeta\in\Gamma(0)}(1-\mod{\zeta})<\infty.
\end{equation}
In this case $\sum_{\gamma\in \Gamma}(1-\mod{\gamma(w)})$ converges for any point $w\in \bb{D}$. An elementary calculation shows that the sum in~\eqref{blaschkesum} converges if and only if the Poincar\'e series $\sum_{\gamma\in\Gamma} \mod{\gamma'(w)}$ converges for any point $w\in \bb{D}$. Automorphism groups for which either, and hence both, of these series converges are said to be of \textit{convergence type}; see Tsuji~\cite[Theorem~XI.3]{tsuji}. It is straightforward to see that if $\Gamma$ is a group of convergence type, then $\Gamma$ is Fuchsian.  

If the series in~\eqref{blaschkesum} converges, then we can form the Blaschke product with zero set $\Gamma(0)$ and this Blaschke product is character automorphic. In Section~\ref{factorization}, we will use these Blaschke products to construct orthonormal bases for the character spaces $H^2_\sigma$. In Section~\ref{forelli} we review some of the results from Forelli~\cite{forelli}. In particular, we revisit the defect space $N$. We assume, as did Forelli, that the space $N$ is finite-dimensional. We recall some duality results for the Hardy spaces $H^p_\Gamma$. These duality results are a generalization of Forelli's original result~\cite[Lemma~3]{forelli}, and the result of Earle and Marden~\cite[Proposition~5]{earlemarden}. The proof we provide is purely functional analytic. We call a group $\Gamma$ \textit{admissible}  if and only if $\Gamma$ is a Fuchsian group of convergence type for which the defect space is finite dimensional. In Section~\ref{rkhs}, we look more closely at the reproducing kernel function for the space $H^2_\sigma$ and show how it is related to the Szeg\"o kernel. These results imply that $H^\infty_\sigma$ is dense in the space $H^p_\sigma$. In section~\ref{hinftyg}, we examine more closely the structure of the operator algebra $H^\infty_\Gamma\subset B(H^2_\Gamma)$. The main results of our paper are proved in Section~\ref{abrahamse}. In particular, we prove a distance formula, Theorem~\ref{distform}, which implies the following analogue of the Nehari theorem:
\begin{thm}
Let $\Gamma$ be Fuchsian group of convergence type with finite dimensional defect space $N$, that is, let $\Gamma$ be admissible. Let $f\in L^\infty_\Gamma$. The distance of $f$ from $H^\infty_\Gamma$ is given by 
\begin{equation*}
\norm{f+H^\infty_\Gamma}=\sup_{\sigma\in \hat{\Gamma}}\norm{(I-P_{H^2_\sigma})M_fP_{H^2_\sigma}}. 
\end{equation*}
\end{thm} 
By specializing to the case where $f\in H^\infty_\Gamma$, this distance formula leads to a generalization of Abrahamse's theorem.
\begin{thm}
Let $\Gamma$ be admissible. Let $z_1,\ldots,z_n\in \bb{D}$ and $w_1,\ldots,w_n\in \bb{C}$. There exists a function $f\in H^\infty_\Gamma$ with $\norm{f}_\infty\leq 1$ such that $f(z_j)=w_j$ if and only if 
\begin{equation*}
A_\sigma :=[(1-w_i\overline{w_j})K^\sigma(z_i,z_j)]_{i,j=1}^n\geq 0 
\end{equation*}
for all $\sigma\in \hat{\Gamma}$. 
\end{thm}
Finally, in Section~\ref{examples} we illustrate our ideas in the context of the two simplest examples of Fuchsian groups.

\section{Factorization and Blaschke products}\label{factorization}

\subsection{Factorization results}
The first step in our work is to establish the analogues of the two main factorization theorems for Hardy spaces: the inner-outer factorization and the Riesz factorization. The inner-outer factorization appears in one form or another in Abrahamse~\cite[Theorem~1.12]{abrahamse2}, Hasumi~\cite[Lemma~1]{hasumi}, and Voichick and Zalcman~\cite{voichickzalcman}. A version of the Riesz factorization is in Abrahamse~\cite[Lemma~5]{abrahamse1}. As we will make frequent use of these results we state and give a short proof of them.

Let $\Gamma$ be an automorphism group. Recall that an element $f\in L^p$ is called \textit{character automorphic} if and only if there exists a character $\sigma\in \hat{\Gamma}$ such that $f\circ \gamma = \sigma(\gamma) f$ for all $\gamma\in \Gamma$. A function $f\in L^p$ is said to be \index{Modulus automorphic function}\textit{modulus automorphic} if and only if $\mod{f}\in L^P_\Gamma$. The absolute value of a character automorphic element of $L^p$ is modulus automorphic. A function $u\in H^p$ is called \textit{outer} if the closed linear span in $H^p$ of $\{z^nu\,:\, n\geq 0\}$ is all of $H^p$. The following result shows that a modulus automorphic outer function is character automorphic.

\begin{prop}\label{outerchar}
If $u\in H^p$ is an outer function and $\mod{u}\in L^p_\Gamma$, then there exists $\sigma\in \hat{\Gamma}$ such that $u\in H^p_\sigma$. 
\end{prop}

\begin{proof}
If $u$ is outer, then $u\circ \gamma$ is outer for all $\gamma\in \Gamma$, since composition by $\gamma$ is continuous and invertible  on $H^p$. Two outer function $u$ and $v$ have equal modulus if and only if there exists a scalar $\lambda\in \bb{T}$ such that $u=\lambda v$, see~\cite[Corollary~6.23]{douglas}. We have $\mod{u\circ \gamma}=\mod{u}\circ \gamma = \mod{u}$ and so there exists $\chi(\gamma)\in \bb{T}$ such that $u\circ \gamma = \sigma(\gamma) u$. We need to show that $\sigma\in \hat{\Gamma}$. If $\gamma_1,\gamma_2\in \Gamma$, then 
\[
\sigma(\gamma_1\gamma_2)u = u\circ (\gamma_1\gamma_2) = (u\circ \gamma_1)\circ \gamma_2
=\sigma(\gamma_1)(u\circ \gamma_2) = \sigma(\gamma_1)\sigma(\gamma_2) u.
\]
Since $u$ is nonzero, $\sigma$ is a character.
\end{proof}

We can now state the two factorization theorems that we need.

\begin{prop}\label{bsufactorization}
Let $\sigma\in\hat{\Gamma}$, let $f\in H^p_\sigma$, and let $f=Bsu$ be the canonical factorization~\cite[Theorem~20 and Chapter~4.4]{helson} of $f$ into a Blaschke product $B$, a singular inner function $s$, and an outer function $u$. Then there exist characters $\sigma_1,\sigma_2,\sigma_3\in \hat{\Gamma}$ such that $\sigma = \sigma_1+\sigma_2+\sigma_3$, $B\in H^\infty_{\sigma_1}$, $u\in H^1_{\sigma_2}$ and $s\in H^\infty_{\sigma_3}$.
\end{prop}

\begin{proof}
Note that if $f(z) = 0$, then $f(\gamma(z)) = 0$ and so the zeros of $f$ are made up of the union of countably many disjoint orbits. The Blaschke product $B$ vanishes precisely on the zero set of $f$. Since $\gamma$ permutes the orbit of a point, we see that $B\circ \gamma$ also vanishes on the zero set of $f$ and so $B\circ \gamma = B C$ where $C$ is inner. A similar argument shows that $B\circ \gamma^{-1} = BD$ with $D$ an inner function. We have,
\begin{align*}
B &= B\circ \gamma\circ \gamma^{-1} = (BC)\circ \gamma^{-1} = (B\circ \gamma^{-1})(C\circ \gamma^{-1}) = BD(C\circ \gamma^{-1}).
\end{align*}
Since $H^\infty$ has no zero divisors we get that $D(C\circ \gamma^{-1}) = 1$ and so $\cl{D} = C\circ \gamma^{-1}\in H^\infty$. This shows that $D$ and $C$ are constant and so $B\circ \gamma = \sigma_1(\gamma) B$ for some scalar $\sigma_1(\gamma)\in \bb{T}$. We now proceed as in the proof of Proposition~\ref{outerchar} to check that $\sigma_1\in \hat{\Gamma}$. 

Since $\mod{f}=\mod{u}\in L^1_\Gamma$, we see by Proposition~\ref{outerchar} that $u\in H^1_{\sigma_3}$ for some character $\sigma_3$. It now follows from the uniqueness of the factorization that $s \in H^\infty_{\sigma_2}$, where $\sigma_2 = \sigma-(\sigma_1+\sigma_3)$.
\end{proof}

One consequence of the proof of Proposition~\ref{bsufactorization} is the following: if the zero set of the Blaschke product $B$ is the orbit $\Gamma(z)$, then $B\in H^\infty_\sigma$ for some $\sigma\in \hat{\Gamma}$. This fact has been observed many times in the literature~\cite{abrahamse1} and we will make this formal later on in this section. 

\begin{prop}[Riesz factorization]\label{riesz2}
Let $\Gamma$ be an automorphism group and let $\sigma\in \hat{\Gamma}$. If $f\in H^1_\sigma$, then there exists  characters $\sigma_1,\sigma_2\in \hat{\Gamma}$ with $\sigma=\sigma_1+2\sigma_2$, an inner function $\phi\in H^\infty_{\sigma_1}$, and an outer function $u\in H^2_{\sigma_2}$   such that $f=\phi u^2$.
\end{prop}

\begin{proof}
It is well known \cite[Theorem~19]{helson} that $f$ has a factorization of the form $\phi u^2$, where $\phi$ is inner and $u\in H^2$ is outer. Note that $\mod{f}^{1/2}=\mod{u}\in L^2_\Gamma$ and so by Proposition~\ref{outerchar}, $u\in H^2_{\sigma_2}$ for some $\sigma_2\in \hat{\Gamma}$. It follows easily that $\phi \in H^\infty_{\sigma_1}$ with $\sigma_1 = \sigma-2\sigma_2$.
\end{proof}

\subsection{Character automorphic Blaschke products}
If $\Gamma$ is of convergence type, then there is a natural class of character-automorphic Blaschke products. These Blaschke products occur naturally when we carry out the inner-outer factorization for elements of $H^1_\Gamma$. For us, the primary use of these Blaschke products is in constructing an orthonormal basis for $H^2_\sigma$. It is easy to check that a Blaschke product $B$ is character automorphic if and only if its zero set is an at most countable union of orbits. We begin by analyzing the convergence of the Blaschke sum~\eqref{blaschkesum}.

\begin{prop}[{\cite[Theorem 2.3.5]{katok}}]\label{stabilizer}
Let $\Gamma$ be a group of automorphisms. If $\fix(\Gamma,0)$ is finite, then $\fix(\Gamma,0)$ is cyclic.
\end{prop}

\begin{prop}\label{blaschkesum2}
Let $\Gamma$ be a group of automorphisms of the disk. Assume that the algebra $H^\infty_\Gamma$ is nontrivial, that is, contains a nonconstant function. The following are true:
\begin{enumerate}
\item For every $w\in \bb{D}$ the stabilizer $\fix(\Gamma,w)$ is a finite,  cyclic group.  
\item For every $w\in \bb{D}$, the Blaschke sum $\sum_{\gamma \in \Gamma}(1-\mod{\gamma(w)})$ converges.
\end{enumerate}
\end{prop}
\begin{proof} Let $w\in \bb{D}$. Choose an automorphism $\phi$ that maps $w$ to $0$.  Let $\Gamma'=\phi \Gamma\phi^{-1}$ and note that $\fix(\Gamma',0) = \fix(\Gamma,0)$ and $\Gamma'(0) = \Gamma(w)$. Also note that the spaces $H^\infty_{\Gamma'}$ and $H^\infty_\Gamma$ are isometrically isomorphic. Therefore, it is enough to consider the case $w=0$.
\begin{enumerate}
\item Every element $\gamma\in \fix(\Gamma,0)$ fixes the origin. Hence, $\gamma$ is a rotation of the disk and there exists a constant $\lambda\in \bb{T}$ such that $\gamma(z)=\lambda z$. Let $f$ be a nonconstant function in $H^\infty_\Gamma$ and let $a_k$, $k\not=0$, be a nonzero Fourier coefficient of $f$. For every $\gamma\in\fix(\Gamma,0)$ we have $f(\gamma(z))=f(z)$ and so $\lambda^ka_k = a_k$. This yields, $\lambda^k=1$ and so $\fix(\Gamma,0)$ is finite. The claim about $\fix(\Gamma,0)$ being cyclic follows from Proposition~\ref{stabilizer}.
\item If $f\in H^\infty_\Gamma$ is nonconstant, then by subtracting a constant we may assume that $f$ is nonzero and vanishes at $0$. Since $f\circ \gamma=f$ we see that $f(\gamma(0))=f(0)=0$ for all $\gamma\in \Gamma$ and so $\sum_{\zeta \in \Gamma(0)}(1-\mod{\zeta})<\infty$. If $\alpha\in \Gamma$, then the cardinality of $\fix(\Gamma,0)$ is equal to the cardinality of $\fix(\Gamma,\alpha(0))$. In fact, the two stabilizers are isomorphic via the map $\gamma\mapsto \alpha\gamma \alpha^{-1}$. Therefore, 
\[
\sum_{\gamma\in \Gamma}(1-\mod{\gamma(0)}) = \mod{\fix(\Gamma,0)}\sum_{\zeta\in \Gamma(0)}(1-\mod{\zeta})<\infty.
\]
\end{enumerate}
\end{proof}

Note that if the series $\sum_{\gamma\in\Gamma}(1-\mod{\gamma(0)})$ does converge, then the argument above shows that $\fix(\Gamma,0)$ is finite and $\Gamma$ is discrete.  

If the stabilizer $\fix(\Gamma,w)$ is finite, and the series $\sum_{\gamma\in\Gamma}(1-\mod{\gamma(w)})$ converges, then we define two Blaschke products $B_{w}$ and $B_{\Gamma,w}$ that arise naturally. First consider the Blaschke product $B_w$ whose zero set is $\Gamma(w)$. We call $B_w$ the \textit{Blaschke product for the orbit $\Gamma(w)$}. Let $m=\mod{\fix(\Gamma,w)}$ and define $B_{\Gamma,w}:=B_w^m$. We call $B_\Gamma$ the \textit{Blaschke product associated to the pair} ($\Gamma$,0). If it is the case that only the identity map fixes $w$, then $m=1$ and $B_{\Gamma,w}=B_w$. If $\gamma\in \Gamma$, then $B_w\circ \gamma$ is a Blaschke product whose zero set is the same as the zero set of $B_w$. It follows, just as in Proposition~\ref{bsufactorization}, that $B_w\circ \gamma = \sigma_w(\gamma) B_w$ for some character $\sigma_w$. We call $\sigma_w$ the \textit{character associated to the pair $(\Gamma,w)$}. We will most often be interested in the case where $w=0$, in which case we denote $B_{\Gamma,0}$ by $B_\Gamma$.

We now construct an orthonormal basis for $H^2_\sigma$ from the Blaschke product $B_0$ and the kernel functions for the spaces $H^2_\sigma$. Our purpose in the next few lemmas is to show that it is enough to construct an orthonormal basis in the case where $H^2_\sigma$ contains at least one function that does not  vanish at the origin. 

\begin{lemma}\label{blaschkecharacter}
Let $\gamma\in \Gamma$ be a generator of the cyclic group $\fix(\Gamma,0)$ and let $\gamma(z)=\lambda z$, where $\lambda\in\bb{T}$. Let $\gamma_n(0)$ be an enumeration of the distinct nonzero elements of $\Gamma(0)$. Then $B_0\circ \gamma  = \lambda B_0$ and $B_\Gamma \circ \gamma = B_\Gamma$. 
\end{lemma}
\begin{proof}
The Blaschke product $B_0$ can be written 
\[
B_0(z) = z\prod_{n=1}^N \frac{\mod{\gamma_n(0)}}{\gamma_n(0)}\frac{\gamma_n(0)-z}{1-\cl{\gamma_n(0)}z}.
\]
Hence, 
\begin{align*}
B_0(\lambda z) &= (\lambda z)\prod_{n=1}^N \frac{\mod{\gamma_n(0)}}{\gamma_n(0)}\frac{\gamma_n(0)-\lambda z}{1-\cl{\gamma_n(0)}\lambda z}= \lambda z\prod_{n=1}^N \frac{\mod{\cl{\lambda}\gamma_n(0)}}{\cl{\lambda}\gamma_n(0)}\frac{\cl{\lambda}\gamma_n(0)-z}{1-\cl{\gamma_n(0)}\lambda z}.
\end{align*}
We have, $\cl{\lambda}\gamma_n(0) = \gamma^{-1}(\gamma_n(0))$ and so $\cl{\lambda}\gamma_n(0)$ is another enumeration of the nonzero points in the orbit $\Gamma(0)$. Hence, $B_0(\lambda z) = \lambda B_0(z)$. If $\mod{\fix(\Gamma,0)}=m$, then $\gamma(z) = e^{\frac{2\pi i}{m}}(z)$ is a generator for $\fix(\Gamma,0)$ and so $B_\Gamma(e^{\frac{2\pi i}{m}}z)=(e^{\frac{2\pi i}{m}}B_0(z))^m = B_\Gamma$ for all $z\in \bb{D}$. 
\end{proof}

\begin{lemma}\label{blaschkedivision}
Let $f\in H^1_\sigma$ be a character automorphic function. Let $m=\mod{\fix(\Gamma,0)}$ and let  $\rho(z) = e^{2\pi i /m} z$ be the generator of $\fix(\Gamma,0)$. Then $\sigma(\rho)=e^{2\pi i r/m}$, $0\leq r\leq m-1$. If $r\not =0$, then $f(0) = 0$. If $r=0$ and $f(0)=0$, then $B_\Gamma|f$.
\end{lemma}

\begin{proof}
Let $\Gamma_1 = \cycgroup{\rho}$ and let $\omega = e^{2\pi i/m}$. Since $\rho$ has order $m$, we see that there exists an integer $r$, $0\leq r\leq m-1$, such that $\sigma(\rho) = \omega^r$. Since $f\circ \rho = \sigma(\rho) f = \omega^r f$, we see that the Fourier expansion of $f$ has the form $f = \sum_{j=0}^\infty a_{mj+r}z^{mj+r}$. If $r\not =0$, then $f(0)=0$. On the other hand if $r=0$, then $z^m|(f-f(0))$. Hence, $B_\Gamma|(f-f(0))$.
\end{proof}

The Hardy space $H^2$ is a reproducing kernel Hilbert space, and the kernel function for $H^2$ is the Szeg\"o kernel 
\[K^S(z,w) = \frac{1}{1-\cl{w}z}.\] 
If $\cal{M}$ is a subspace of $H^2$, then $\cal{M}$ is also a reproducing kernel Hilbert space and we denote its kernel function $K^\cal{M}$. We denote the kernel function for $H^2_\sigma$ by $K^\sigma$ and the kernel function for $H^2_\Gamma$ by $K^\Gamma$. Given a point $w\in\bb{D}$ and a kernel function $K$ we denote by $k_w$, the kernel function at the point $w$, that is, $k_w$ is the unique element of $\cal{M}$ such that $\inp{f}{k_w} = f(w)$ for all $f\in \cal{M}$. It follows that $\norm{k_w}^2 = K(w,w)$. If $k_w\not =0$, then we define the normalized kernel function by $\kappa_w = K(w,w)^{-1/2}k_w$. Otherwise, we set $\kappa_w = 0$.  Recall that $\Gamma$ is of convergence type if the sum in~\eqref{blaschkesum} converges.

\begin{prop}\label{onb}
Let $\Gamma$ be a Fuchsian group of convergence type. Let $m=\mod{\fix(\Gamma,0)}$. Assume that $\sigma\in \hat{\Gamma}$ and that there exists a function $f\in H^2_\sigma$ with $f(0)\not = 0$. For $n\geq 0$, let $\sigma_n = \sigma-mn\sigma_0$,  let $H^2_n = H^2_{\sigma_n}$ and let $\kappa_0^{(n)}$ denote the normalized kernel function for $H^2_n$. Then the nonzero elements of the set 
\[\cal{E}_\sigma :=\{B_\Gamma^{n}\kappa^{(n)}_0\,:\,n\geq 0\},\]
is an orthonormal basis for $H^2_\sigma$.
\end{prop}

\begin{proof}
Since $B_0$ is an inner function and $\kappa^{(n)}_0$ is the normalized kernel function we see that the nonzero elements of $\cal{E}_\sigma$ have norm 1. Since $B_0(0)=0$, we see that the nonzero elements of $\cal{E}_\sigma$ are orthogonal.

If $\rho$ is the generator of $\fix(\Gamma,0)$, then $\sigma(\rho) =1$, by Lemma~\ref{blaschkedivision}. Combining this with Lemma~\ref{blaschkecharacter} we see that $\sigma_n(\rho) = 1$ for $n\geq 1$.

Now let $f\in H^2_\sigma$ with $f\perp \cal{E}_\sigma$. Note that $\kappa_0^{(n)}=0$ implies that every element of $H^2_n$ vanishes at the origin. Since $f\perp \kappa^{(0)}_0$, we see that $f(0)=0$ and so $B_\Gamma|f$, by Lemma~\ref{blaschkedivision}. Therefore, $f=B_\Gamma f_1$ and by composing with elements of $\Gamma$ we see that $f_1\in H^2_{\sigma-m\sigma_0}=H^2_{\sigma_1}$. If $k^{(1)}_0=0$, then every element of $H^2_1$ vanishes at the origin and so $B_\Gamma|f_1$. If $k^{(1)}_0\not=0$, then the fact that $f\perp B_\Gamma\kappa^{(1)}_0$ gives $f_1(0)=0$ and $B_\Gamma|f_1$. Repeating this process we see that $B_\Gamma^{n}$ divides $f$ for all $n\geq 0$. Hence, $f=0$.
\end{proof}

\section{The Forelli projection}\label{forelli}

\subsection{The conditional expectation}
To begin this section we review some ideas from the work of Forelli~\cite{forelli}. Let $\Gamma$ be an automorphism group. Let $\mathfrak{M}$ denote the collection of Lebesgue measurable subsets of $\bb{T}$. For $E\in\mathfrak{M}$, define $E$ to be $\Gamma$-invariant if and only if $m(E\triangle \gamma^{-1}(E))=0$ for every $\gamma\in \Gamma$. Forelli proved the following fact: the collection of $\Gamma$-invariant sets $\mathfrak{M}_\Gamma$ is a sub-sigma-algebra of $\mathfrak{M}$, and the subspace $L^p(\bb{T},\mathfrak{M}_\Gamma,m)$ is equal to $L^p_\Gamma$. 

It follows, from standard arguments in probability theory, that there is the unique conditional expectation $\Psi_p:L^p(\bb{T},\mathfrak{M},m)\to L^p(\bb{T},\mathfrak{M}_\Gamma,m)$. The conditional expectation $\Psi_p$ is $\wk$-continuous for $1<p\leq \infty$, and $\Psi_p|L^q = \Psi_q$ for $p\leq q$. Therefore, we omit the subscript $p$ and denote $\Psi_p$ by $\Psi$. Note that the map $\Psi$ is a projection ($\Psi^2=\Psi$), and is selfadjoint ($\Psi^\ast=\Psi$) when $p=2$. The map $\Psi$ has the property that 
\[\int_E \Psi(f) d\mu = \int_E fd\mu,\]
for $E\in \mathfrak{M}_\Gamma$. Given $f\in L^p(\bb{T},\mathfrak{M},d\mu)$ and $g\in L^q(\bb{T},\mathfrak{M}_\Gamma,d\mu)$ we have
\[\Psi(fg) = \Psi(f)g.\]

In the classical setting, the space $L^2$ decomposes as a direct sum of the Hardy space $H^2$ and the space $\cl{H^2_0}$, that is, the set of complex conjugates of functions in $H^2$ that have mean value 0. In our setting we must contend with the fact that the direct sum $H^2_\Gamma\oplus \cl{H^2_{0,\Gamma}}$ may not be all of $L^2_\Gamma$. The orthogonal complement of $H^2_\Gamma\oplus \cl{H^2_{0,\Gamma}}$ in $L^2_\Gamma$ is called the \textit{defect space}. If $\Gamma$ is the group of deck transformations that arise from a covering map of a $g$-holed region, then $N$ is a $g$-dimensional subspace. Forelli captures the defect space without reference to multiply connected domains. We now outline these ideas, as they are central to our work. We also feel that our presentation is sufficiently different from~\cite{forelli} to merit inclusion.

Let $k_z\in H^2$ denote the Szeg\"o kernel for the Hardy space. To begin with let $f\in L^2_\Gamma$ and let $f=g+\cl{h}+c$, where $g,h\in H^2_0$ and $c = \int f dm$. While the functions $g,h$ are not $\Gamma$-invariant, they are ``nearly'' so. We have 
\[g+\cl{h}+c = f = f\circ \gamma = g\circ \gamma + \cl{h\circ\gamma} +c.\]
Rearranging this last equation gives
\[g\circ \gamma - g = \cl{h-h\circ \gamma}.\]
Since $\gamma\in \Gamma$ is analytic we see that $g\circ \gamma$ and $h\circ \gamma$ are both in $H^2$ and so $g\circ\gamma - g = \cl{h-h\circ \gamma} = c(f,\gamma)$. Where $c(f,\gamma)$ is a constant that depends on $f$ and $\gamma$. 
The constant is easily computed by integrating 
\[
c(f,\gamma)= \int g\circ \gamma - g  = g(\gamma(0))-g(0)  = \inp{g}{k_{\gamma(0)}}.
\]
Similarly, $c(f,\gamma) = \inp{\cl{h}}{\cl{k_{\gamma(0)}}}$. Combining these two equations, we get 
\begin{align*}
2c(f,\gamma) &= \inp{g+\cl{h}}{k_{\gamma(0)}-\cl{k_{\gamma(0)}}}  \\
&= \inp{\Psi(f)}{k_{\gamma(0)}-\cl{k_{\gamma(0)}}} = \inp{f}{\Psi(k_{\gamma(0)}-\cl{k_{\gamma(0)}})}.
\end{align*}

If we set $v_\gamma = -i\Psi(k_{\gamma(0)}-\cl{k_{\gamma(0)}})$, then $v_\gamma$ is a positive element of $L^\infty_\Gamma$ with the property that 
\[
\int fv_\gamma = 2ic(f,\gamma)
\]
for all $\gamma\in \Gamma$. Note that $-i(k_\gamma(0)-\cl{k_{\gamma(0)}})=P_{\gamma(0)}^\ast$, the conjugate of the Poisson kernel at the point $\gamma(0)$. Following Forelli we define the \index{Defect space}\textit{defect space} $N$ by
\[
N:=\linspan\{v_\gamma\,:\,\gamma\in \Gamma\}.
\]

If $f\in L^2_\Gamma$ is orthogonal to $v_\gamma$, then $g\circ\gamma - g = 0 = h-h\circ\gamma$.  Hence, $g$ and $h\in H^2_{0,\Gamma}$ and we get that $L^2_\Gamma = H^2_\Gamma\oplus \cl{H^2_{0,\Gamma}}\oplus [N]$. 

If $\gamma_1,\gamma_2\in \Gamma$, then 
\begin{align*}
c_{\gamma_1\circ\gamma_2} & = g\circ (\gamma_1\circ\gamma_2) - g =  g\circ (\gamma_1\circ\gamma_2) - g \circ \gamma_2 + g\circ \gamma_2 -  g \\
& = (g\circ \gamma_1 - g)\circ \gamma_2 +(g\circ \gamma_2 - g) = c_{\gamma_1}+c_{\gamma_2}
\end{align*}

It follows that 
\[
\int fv_{\gamma_1\circ\gamma_2} = 2ic_{\gamma_1\circ\gamma_2} = 2i(c_{\gamma_1}+c_{\gamma_2}) = \int f(v_{\gamma_1}+v_{\gamma_2}). 
\]
Since this is true for all $f\in L^2_\Gamma$ we get $v_{\gamma_1\circ\gamma_2}=v_{\gamma_1}+v_{\gamma_2}$. 

What all of this shows is that the map $\gamma \mapsto v_\gamma$ is a homomorphism from $\Gamma$ into the additive group $N$. This homomorphism must factor through the commutator subgroup $[\Gamma,\Gamma]$ to give a homomorphism from $\Gamma/[\Gamma,\Gamma]$ into $N$. 

Given a set of generators $\{\dot{\gamma_s}\,:\,s\in S\}\subseteq  \Gamma/[\Gamma,\Gamma]$ the vectors $v_s := v_{\gamma_s}$, for $s\in S$, span the space $N$. If the group $\Gamma$ is finitely generated, then the space $N$ is finite dimensional and the dimension of $N$ is smaller than the minimal number of generators of $\Gamma/[\Gamma,\Gamma]$. Note that if an element $\gamma\in \Gamma$ has finite order, say $m$, then $mv_\gamma = v_{\gamma^m} = 0$ and so $v_\gamma = 0$. If $\Gamma/[\Gamma,\Gamma]$ is generated by elements of finite order, then $N$ is trivial. 

In the case where $\Gamma$ is the group of deck transformations associated with a universal covering map of a multiply connected domain, Forelli showed that  the dimension of $N$ is equal to the minimal number of generators of $\Gamma/[\Gamma,\Gamma]$. Many results follow from this equality, including a corona theorem for $H^\infty_\Gamma$, the key stepping stone being the construction of a bounded projection $P:H^\infty\to H^\infty_\Gamma$. We do not assume, and do not need, this stronger condition. 

We are now in a position to state what kind of groups we will deal with in this paper. We will call a group $\Gamma$ {\bfseries\itshape admissible} if and only if $\Gamma$ is of convergence type and the defect space $N$ is finite dimensional.

These are natural, nontrivial conditions from a function theory point of view. The assumption about $N$ means that $H^p_\Gamma+N$ is a closed subspace of $L^p$ for all $1\leq p \leq \infty$. It is useful to keep in mind that $N\subset L^\infty_\Gamma\subseteq L^p_\Gamma$ for all $1\leq p\leq \infty$. Note that when $q>1$, the spaces $H^q_\Gamma$ and $H^q_\Gamma+N$ are also \wk{} closed. This is a simple consequence of the fact that \wk{} limits preserve point values and the fact that $N$ is assumed finite dimensional. 

\subsection{Duality and density results}
Since duality arguments will play a central role in our interpolation results we would like to gather some results on the duality between the different $H^p_\Gamma$ spaces. The results contained in the next proposition have appeared in various forms in the literature~\cite[Theorem~11.1]{ahernsarason}, \cite[Proposition~5]{earlemarden}, and \cite[Lemma~3]{forelli}. We provide, mostly for the sake of completeness, an elementary proof that depends only on standard duality arguments and the expectation $\Psi$. 

\begin{prop}\label{duals}
Let $\Gamma$ be an admissible group. The image of $H^p$ under the conditional expectation $\Psi$ is $H^p_\Gamma+N$. For $1\leq p< \infty$, the dual of $L_\Gamma^p$ can be identified with $L^q_\Gamma$, where $q=\frac{p}{p-1}$. In this identification the following are true:
\begin{enumerate}
\item $(H^p_\Gamma)^\perp = H^q_{0,\Gamma}+N$
\item $(H^q_\Gamma)_\perp = H^p_{0,\Gamma}+N$
\item $(H^q_\Gamma+N)_\perp = H^p_{0,\Gamma}$
\item $(H^p_\Gamma+N)^\perp = H^q_{0,\Gamma}$
\end{enumerate}
\end{prop}

\begin{proof}
The statement about $L^p_\Gamma$ spaces follows from standard facts about the $L^p$ spaces of a probability measure. 

For $p=2$, we have seen already that $L^2_\Gamma = H^2_\Gamma\oplus \cl{H^2_{0,\Gamma}}\oplus N$ and the above results are valid. We have
\[
\int \Psi(h)g = \int hg = 0, 
\]
for all $h\in H^2$, $g\in H^2_{0,\Gamma}$, and so $\Psi(H^2) = L^2_\Gamma\ominus \cl{H^2_{0,\Gamma}}=H^2_\Gamma\oplus N$. 

Next consider the case $1\leq p <\infty$. We will show that $H^p_\Gamma = (H^q_{0,\Gamma}+N)_\perp$ and the remaining results will follow either by duality or by a similar argument. Note that the spaces $H^q_\Gamma$ and $H^q+N$ are \wk{} closed in $L^q$.

We first show that $\Psi(H^p)=H^p_\Gamma+N$. Let $f\in L^q_\Gamma$. Note that $f\in\Psi(H^p)^\perp$ if and only if $\int f\Psi(g)=\int fg =0$ for all $g\in H^p$ if and only if $f\in H^q_0\cap L^q_\Gamma =H^q_{0,\Gamma}$. Hence, $\Psi(H^p)^\perp = H^q_{0,\Gamma}$. 

Let $h\in H^p_\Gamma$, $v\in N$ and $g\in H^q_{0,\Gamma}$. Since $v_\gamma\in L^\infty_\Gamma$ and $g$ is analytic we get
\[
\int gv_\gamma = \int g\Psi(k_{\gamma(0)}-\cl{k_{\gamma(0)}}) = \int g(k_{\gamma(0)}-\cl{k_{\gamma(0)}}) = g(\gamma(0))- g(0) = 0.
\]
It follows that
\[
\int g(h+v) = \int gh+\int gv = 0, 
\]
which gives $H^p_\Gamma+N\subseteq (H^q_{0,\Gamma})_\perp$. This also yields
\[
\Psi(H^p) = (\Psi(H^p)^\perp)_\perp = (H^q_{0,\Gamma})_\perp \supseteq H^p_\Gamma+N. 
\]

Let $g\in (H^p_{\Gamma}+N)^\perp$ and let $f\in H^p$. We will show that 
\[
\int gf = \int g\Psi(f) = 0, 
\]
which will establish the fact that $(H^p_{\Gamma}+N)^\perp \subseteq \Psi(H^p)^\perp$. We know that $\int g(h+v) = 0$ for all $h \in H^p_{\Gamma}$ and $v\in N$. Let $f_n\in H^\infty$ and suppose that $f_n\to f$ in the $L^p$ norm. Since $f_n\in H^2$, we can write $\Psi(f_n) = h_n+v_n$, where $h_n\in H^2_\Gamma$ and $v_n\in N$. However, $v_n\in L^\infty_\Gamma$ and so $h_n\in H^\infty_\Gamma\subseteq H^p_\Gamma$. Therefore, 
\[
\int g(h_n+v_n) =\int g\Psi(f_n) = 0
\]
for all $n$, since $g\in \Psi(H^p)^\perp$. Hence, 
\[
\int g\Psi(f) = \lim_{n\to\infty}\int g\Psi(f_n) = 0.
\]

We have established that $\Psi(H^p)^\perp = H^q_{0,\Gamma}$ and that $\Psi(H^p)^\perp \supseteq (H^p_\Gamma+N)^\perp$. This combined with the fact that $\Psi(H^p) \supseteq H^p_\Gamma+N$ shows us that $\Psi(H^p) = H^p_\Gamma+N$ and $(H^p_\Gamma+N)^\perp = H^q_{0,\Gamma}$.
\end{proof}

\begin{prop}\label{hinftygclosure}
For $1\leq p< \infty$, $[H^\infty_\Gamma]_p=H^p_\Gamma$. 
\end{prop}

\begin{proof}
Consider the case $p\geq 2$. Let $f\in H^p_\Gamma$ and let $f_n\in H^\infty$ converge to $f$ in $L^p$. Since the $L^p$ norm dominates the $L^2$ norm we see that $f_n\to f$ in $L^2$. If we project $f_n$ onto $\Psi(f_n)$, then there exists $g_n\in H^\infty_\Gamma$ and $v_n\in N$ such that $\Psi(f_n)=g_n+v_n$ and $g_n+v_n\to f$ in the $L^p$ norm and also in the $L^2$ norm. In particular, $\norm{v_n}_2 \to 0$. However, on the finite dimensional space $N$, the $L^2$ norm and the $L^p$ norm are equivalent and so $\norm{v_n}_p \to 0$. It follows that $g_n\to f$ in $L^p$.

Let $1\leq p<2$, view $H^\infty_\Gamma$ as a subspace of $H^p_\Gamma$, and let $(H^\infty_\Gamma)^\perp$ denote the annihilator of $H^\infty_\Gamma$ in $L^q_\Gamma$, where $2\leq q$ and $q^{-1}+p^{-1}=1$. We have
\begin{align*}
(H^\infty_\Gamma)^\perp &=\bigg\{f\in L^q_\Gamma\,:\, \int fg=0 \text{ for all } g\in H^\infty_\Gamma\bigg\}\\
&\subseteq \left\{f\in L^2_\Gamma\,:\, \int fg=0 \text{ for all } g\in H^\infty_\Gamma\right\}\\
&=L^2_\Gamma\ominus \left[\cl{H^\infty_\Gamma}\right]=L^2_\Gamma\ominus \cl{H^2_\Gamma} =H^2_{0,\Gamma}\oplus N.
\end{align*}
Hence, $(H^\infty_\Gamma)^\perp\subseteq L^q_\Gamma\cap(H^2_{0,\Gamma}+N) = H^q_{0,\Gamma}+N$, since $N\subseteq L^\infty_\Gamma\subseteq L^q_\Gamma$. The reverse inclusion follows from the inclusions $H^q_{0,\Gamma}+N\subseteq H^2_{0,\Gamma}\oplus N$ and $H^\infty_\Gamma\subseteq H^2_\Gamma$, and the result for $p=2$. We have shown that $(H^\infty_\Gamma)^\perp = H^q_{0,\Gamma}+N$. It follows from Proposition~\ref{duals} that $[H^\infty_\Gamma]_p =((H^\infty_\Gamma)^\perp)_\perp = (H^q_{0,\Gamma}+N)_\perp = H^p_\Gamma$.
\end{proof}

We now prove a density result that is central to the proof of our distance formula and interpolation theorem. Recall that a function $u\in H^2$ is called \textit{outer} if and only if $[H^\infty u]=H^2$. If $\sigma\in \hat{\Gamma}$ is a character and $u\in H^2_\sigma$ is an outer function, then $[H^\infty_\Gamma u]\subseteq H^2_\sigma$. The purpose of the next two results is to prove the reverse containment. Let us temporarily adopt the following notation: if $\norm{u}_2=1$, then we call $u$ normalized; and let the $L^2$ norm induced by the measure $\mod{u}^2 dm$ be denoted $\norm{}_{2,u}$.

\begin{lemma}\label{directsum}
Let $\Gamma$ be an admissible group, let $u\in H^2$ be an outer function and let $N$ be finite dimensional. Then $[H^\infty_\Gamma u]\cap [Nu] =\{0\}$.
\end{lemma}

\begin{proof}
We can assume that $u$ is normalized. Let $f \in [H^\infty_\Gamma u]\cap [Nu]$ and note that $f\in H^2$. Since $u$ is an outer function, there exists a sequence $f_n\in H^\infty_\Gamma$ and $v_n\in N$ such that 
$\lim_{n\to\infty}\norm{f-f_nu}_2 = \lim_{n\to\infty}\norm{f-v_nu}_2=0$.  The space $N$ is finite dimensional and so $\norm{}_2$ is equivalent to $\norm{}_{2,u}$. Hence, the sequence $v_n$ converges, say to $v$, in $L^2$. There exists a sequence $h_n\in H^\infty$ such that $\lim_{n\to\infty}\norm{h_n u -1}_2=0$. We have, 
$\lim_{n\to\infty} fh_n = \lim_{n\to\infty} v_nu h_n = \lim_{n\to\infty} v_n \in H^2_\Gamma$. This implies $v=0$ and so $f=0$.
\end{proof}

We are now in a position to prove our main lemma about cyclic subspaces. 

\begin{thm}\label{cyclic}
Let $\Gamma$ be an admissible group, let $\sigma\in \hat{\Gamma}$ and let $u\in H^2_\sigma$ be an outer function. The cyclic subspace generated by $H^\infty_\Gamma$ and $u$ is equal to $H^2_\sigma$. 
\end{thm}

\begin{proof}
One inclusion is straightforward, namely $[H^\infty_\Gamma u] \subseteq H^2_\sigma$. 

For the converse, let $f\in H^2_\sigma \ominus [H^\infty_\Gamma u]$. We will show that $f=0$. Since $u$ is an outer function, there exists a sequence of functions $f_n\in H^\infty$ such that $\lim_{n\to\infty}\norm{f_nu-f}_2=0$. Let $\Psi$ be the conditional expectation from $L^1$ onto $L^1_\Gamma$. Note that $\mod{u}^2,f\cl{u}$ and $\mod{f}^2\in L^1_\Gamma$. Using the fact that $\int \Psi(f)g = \int fg$ for $f\in L^\infty$ and $g\in L^1_\Gamma$ we get
\begin{align*}
\norm{\Psi(f_n)u-f}_2^2&=\int \mod{\Psi(f_n)}^2\mod{u}^2-\Psi(f_n)u\cl{f} - \cl{\Psi(f_n)u}f+\mod{f}^2\,dm\\
&= \int \cl{f_n}\Psi(f_n)\mod{u}^2-\Psi(f_n)u\cl{f}-\cl{f_n u}f+\mod{f}^2\,dm\\
& = \int (\Psi(f_n)u-f)\cl{(f_nu-f)} \,dm = \inp{\Psi(f_n)u-f}{f_nu-f} \\
&\leq \norm{\Psi(f_n)u-f}_2\norm{f_nu-f}_2.
\end{align*}
Hence, $\lim_{n\to\infty}\norm{\Psi(f_n)u-f}_2\leq \lim_{n\to\infty}\norm{f_nu-f}_2=0$.

By Proposition~\ref{duals} we can write $\Psi(f_n) = g_n+v_n$, where $g_n\in H^\infty_\Gamma$ and $v_n\in N$. Hence, $f\in [H^\infty_\Gamma u]+[Nu]$. Since $N$ is finite dimensional, by Lemma~\ref{directsum}, we can define an equivalent norm on $[H^\infty_\Gamma u]+[Nu]$ by $\norm{h+w}'=\norm{h}_2+\norm{w}_2$. Since $f\perp g_nu$, it follows that $\lim_{n\to\infty}\norm{v_nu -f}'=0$, and $f\in [Nu]\cap [H^\infty_\Gamma u]=\{0\}$. 
\end{proof}

\section{Reproducing kernel functions for character spaces}\label{rkhs}

The purpose of the results in this section is to establish connections between the reproducing kernel $K^\sigma$ and the Szeg\"o kernel  $K^S$. As a corollary to these results we will see that $H^\infty_\sigma$ is dense in $H^1_\sigma$. 

\begin{lemma}\label{intersection}
If $\Gamma$ is an admissible group, then $(H^1_\Gamma+N)\cap(\cl{H^1_\Gamma}+N) = N+\bb{C}$.
\end{lemma}

\begin{proof}
Let $f\in (H^1_\Gamma+N)\cap(\cl{H^1_\Gamma}+N)$. We can write $f = g+v=\cl{h}+w$, where $g,h\in H^1_{0,\Gamma}$ and $v,w\in N$. We have $g-\cl{h} = w-v\in N$. By Proposition~\ref{hinftygclosure} we can choose sequences $g_n,h_n\in H^\infty_\Gamma$ such that $g_n\to g$, $h_n\to h$ in $L^1_\Gamma$. Note that $g_n-\cl{h_n}\in H^2_{\Gamma}\oplus\cl{H^2_\Gamma}$ and so $g_n-\cl{h_n}\perp N$. Therefore, $\norm{g-\cl{h}}_1= 0$ which yields $g\in H^1_\Gamma\cap \cl{H^1_\Gamma}$. Hence, $g$ is constant and $f\in N+\bb{C}$.
\end{proof}

\begin{prop}\label{boundedkernel1}
Let $\Gamma$ be an admissible group, let $\sigma\in\hat{\Gamma}$, let $\cal{N}:= H^2_\sigma\ominus [H^\infty_{0,\Gamma}H^2_\sigma]$. If $f\in\cal{N}$, then $\mod{f}^2\in N+\bb{C}$. In particular, the kernel function at the origin $k_0^\sigma$ is bounded. 
\end{prop}

\begin{proof}
Let $f\in \cal{N}$ and let $h\in H^\infty_{0,\Gamma}$. We have $f\perp hf$ 
\[\int \mod{f}^2h = 0,\]
for all $h\in H^\infty_{0,\Gamma}$. By taking the complex conjugate we also get
\[\int \mod{f}^2\cl{h} = 0,\]
and so $\mod{f}^2\in N+\bb{C}$, by Lemma~\ref{intersection}. Since $N\subseteq L^\infty_\Gamma$, this implies that $\mod{f}$ is bounded. 

To prove the claim about the kernel function, note that $k_0^\sigma \perp [H^\infty_{0,\Gamma}H^2_\sigma]$. 
\end{proof}

\begin{cor}\label{density}
Let $\Gamma$ be an admissible group. If $\sigma\in \hat{\Gamma}$, then $H^\infty_\sigma $ is dense in $H^p_\sigma$.
\end{cor}

\begin{proof}
It is enough to prove the result in the case $p=1$. Recall that, in Proposition~\ref{onb}, we had constructed an orthonormal basis for $H^2_\sigma$. Each of the elements of this basis was the product of a Blaschke product and a normalized kernel function at the origin for some character space. These are bounded by Proposition~\ref{boundedkernel1} and we see that $H^\infty_\sigma$ is dense in $H^2_\sigma$. If $f\in H^1_\sigma$, then by Riesz factorization (see Proposition~\ref{riesz2}) there exists $g\in H^2_{\sigma_1}$ and $h\in H^2_{\sigma-\sigma_1}$ such that $f=gh$. Choose a sequence $g_n\in H^\infty_{\sigma_1}$ and $h_n\in H^\infty_{\sigma-\sigma_1}$ such that $\norm{g_n-g}_2\to 0$ and $\norm{h_n-h}_2\to 0$ in $H^2$. It follows that $f$ is the limit of $g_nh_n\in H^\infty_\sigma$ and so $H^\infty_\sigma$ is dense in $H^1_\sigma$. 
\end{proof}

We now indicate how the kernel function for $H^2_\sigma$ is related to the Szeg\"o kernel. We recall the notion of a positive semidefinite function. Given a set $X$ a function $K:X\times X\to \bb{C}$ is called \textit{positive semidefinite} if and only if for every finite set of points $\{x_1,\ldots,x_n\}\subseteq X$, and finite set of scalars $\{\alpha_1,\ldots, \alpha_n\}\subseteq \bb{C}$ we have $\sum_{i,j=1	}^n \cl{\alpha_i}\alpha_j K(x_i,x_j)\geq 0$.

\begin{cor}\label{boundedkernel}
Let $\Gamma$ be an admissible group, let $w\in \bb{D}$ and let $B_\Gamma$ be the Blaschke product associated to the pair $(\Gamma,0)$. Then there exists a constant $C$, which depends only on the group $\Gamma$, such that the following are true:
\begin{enumerate}
\item\label{ineq1} $\norm{k_w^\sigma}_\infty\leq C^2(1-\mod{B_{\Gamma}(w)})^{-1}$ for every $\sigma\in \hat{\Gamma}$. 
\item $C^2K^S(B_\Gamma(z),B_\Gamma(w))-K^\sigma(z,w)$ is a positive semidefinite function on $\bb{D}\times \bb{D}$.
\end{enumerate}
\end{cor}
\begin{proof}
Since $N+\bb{C}$ is finite dimensional we know that there is a constant $C$ such that $\norm{f}_\infty\leq C \norm{f}_2$ for all $f\in N+\bb{C}$. Hence, $\norm{k_0^\sigma}_\infty\leq C\norm{k_0^\sigma}_2$ and $\norm{\kappa_0^\sigma}_\infty \leq C$. Note that this bound is independent of the character $\sigma$.  

From Proposition~\ref{onb} we have an orthonormal basis for the space $H^2_\sigma$ consisting of the nonzero elements of the set $\cal{E}_\sigma=\{B_\Gamma^n\kappa_0^{(n)}\,:\, n\geq 0\}$. For ease of notation let $B=B_\Gamma$ and $f_n = \kappa^{(n)}_0$. We have
\[
K^\sigma(z,w)  = \sum_{n=0}^\infty B^{n}(z)f_n(z)\cl{B^{n}(w)f_n(w)}
\]
We now prove both claims made above. 
\begin{enumerate}
\item Taking the absolute value and using the fact that $\norm{\kappa_0^\sigma}_\infty \leq C$ we get
\begin{equation*}
\mod{K^\sigma(z,w)}\leq C^2\sum_{n=0}^\infty\mod{B(z)}^{n}\mod{B(w)}^{n} \leq  \frac{C^2}{1-\mod{B(z)B(w)}}. 
\end{equation*}
Therefore, $\mod{k_w^\sigma(z)}\leq C^2(1-\mod{B(z)B(w)})^{-1}\leq C^2(1-\mod{B(w)})^{-1}$. Taking the supremum over $z\in \bb{D}$ yields the inequality in~\ref{ineq1}.
\item An element $f\in H^\infty$ has norm at most $C$ if and only if $C^2-f(z)\cl{f(w)}$ is a positive semidefinite function, this follows from the fact that $H^\infty$ is the multiplier algebra of $H^2$. Let $\alpha_1,\ldots,\alpha_m\in \bb{C}$ and $z_1,\ldots,z_m\in \bb{D}$. 
Now
\begin{align*}
&\sum_{i,j=1}^m\cl{\alpha_i}\alpha_j(C^2K^S(B(z_i),B(z_j))-K^\sigma(z_i,z_j))\\
&= \sum_{i,j=1}^m\sum_{n=0}^\infty \cl{\alpha_i}\alpha_j(C^2B(z_i)^{n}\cl{B(z_j)}^{n}-f_n(z_i)\cl{f_n(z_j)}B(z_i)^{n}\cl{B(z_j)}^{n})\\
&=\sum_{n=0}^\infty\sum_{i,j=1}^m \cl{\alpha_i}\alpha_j(C^2-f_n(z_i)\cl{f_n(z_j)})B(z_i)^{n}\cl{B(z_j)}^{n} 
\end{align*}
which is non-negative, since $\norm{f_n}_\infty \leq C$. 
\end{enumerate}
\end{proof}

\section{$H^\infty_\Gamma$ as an operator algebra}\label{hinftyg}

\subsection{Commutants and multiplier algebras}
In this section we provide some basic information about the operator algebra structure of $H^\infty_\Gamma$ when viewed as the multiplier algebra of $H^2_\sigma$. 

\begin{prop}\label{multalg}
If $\sigma\in \hat{\Gamma}$, then the multiplier algebra of $H^2_\sigma$ is $H^\infty_\Gamma$.  
\end{prop}

\begin{proof}
If $f\in H^\infty_\Gamma$ and $g\in H^2_\sigma$, then it is easily checked that $fg\in H^2_\sigma$. 

For the converse we can assume that no inner function divides $H^2_\sigma$. Hence, $k_w^\sigma\not =0$ for all $w\in \bb{D}$ and we see that a multiplier $f$ of $H^2_\sigma$ is a bounded function. By composing with elements of the group and comparing characters we see that $f\in H^\infty_\Gamma$.

It is a well-known fact~\cite[Theorem~6.3]{paulsen} that $\norm{f}_{\mult}\geq \norm{f}_\infty$. On the other hand, given $f\in H^\infty_\Gamma$ and $h\in H^2_\sigma$, we have
\[\int \mod{fh}^2\,dm\leq \norm{f}_\infty^2\norm{h}_2^2.\]
Hence, $\norm{f}_{\mult} = \norm{f}_\infty$.
\end{proof}

We recall the notions of \textit{reflexivity} and \textit{hyperreflexivity} for an operator algebra $\cal{A}\subseteq B(\cal{H})$. An algebra $\cal{A}\subseteq B(\cal{H})$ is called \textit{reflexive} if and only if $\alg(\lat(\cal{A})) = \cal{A}$. A stronger property is that of hyper-reflexivity~\cite{arveson}. The algebra $\cal{A}$ is \textit{hyper-reflexive} if and only if there exists a constant $C$ such that $\norm{T+\cal{A}}\leq C \sup_{P\in \lat(\cal{A})}\norm{(I-P)TP}$ for all $T\in B(\cal{H})$.

If $\cal{L}$ is the lattice in $B(H^2)$ of subspaces of the form $\phi H^2$, where $\phi$ ranges over the set of inner functions, then $\alg(\lat(\cal{L}))=H^\infty$ and so $H^\infty$ is reflexive. We now prove a more general fact that implies that $H^\infty_\Gamma$ is reflexive as a subalgebra of $H^2_\Gamma$.

\begin{prop}
Let $\cal{H}$ be a reproducing kernel Hilbert space. If $\cal{A}$ denotes the multiplier algebra of $\cal{H}$, then $\cal{A}$ is a reflexive operator algebra in $B(\cal{H})$. 
\end{prop}

\begin{proof}
Suppose that $T\in B(\cal{H})$ and $T(\cal{M})\subseteq \cal{M}$ for all $\cal{M}\in \lat(\cal{A})$. The subspace spanned by the kernel function $k_x$ is invariant for $\cal{A}^\ast$ and so $T^\ast k_x \in \linspan{\{k_x\}}$. Hence, there exists a constant $\phi(x)\in \bb{C}$ such that $T^\ast k_x = \phi(x) k_x$ and consequently $T=M_{\overline{\phi}}$.
\end{proof}

Our next proposition implies that $H^\infty_\Gamma$ is equal to its own commutant when viewed as a subalgebra of $B(H^2_\Gamma)$.

\begin{prop}
Let $\cal{M}$ be a subspace of $H^2$ such that $1\in \cal{M}$
and let $\cal{A}$ be the multiplier algebra of this subspace. If $[\cal{A}]=\cal{M}$, then $\cal{A}'=\cal{A}$, when $\cal{A}$ is represented as multiplication operators on $\cal{M}$.
\end{prop}

\begin{proof}
Let $T\in \cal{A}'$ and set $T(1) = h$. Consider the
action of $T^\ast$ on the kernel function $k_x\in \cal{M}$. Let $g\in \cal{A}$ and
compute 
\begin{align*}
\nonumber\inp{T^\ast k_x}{g} &= \inp{k_x}{TM_g1} = \inp{k_x}{M_gT(1)} = \inp{k_x}{gh} \\
&=\overline{g(x)h(x)}=\inp{\overline{h(x)}k_x}{g}
\end{align*}
and so $T^\ast k_x = \overline{h(x)}k_x$. This forces $h\in \cal{A}$ and $T=M_h$. 
\end{proof}

If $\gamma:\bb{D}\to\bb{D}$ is an automorphism, then an easy calculation,
\[
\inp{f}{C_\gamma^\ast k_z} = \inp{C_\gamma f}{k_z} = f(\gamma(z)) = \inp{f}{k_{\gamma(z)}},
\]
shows that $C_\gamma^\ast k_z = k_{\gamma(z)}$. Given a Fuchsian group $\Gamma$, let us denote by $\alg(S,\Gamma)$ the smallest subalgebra of $B(H^2)$ that contains the unilateral shift $S$ and the group of composition operators $\{C_\gamma\,:\,\gamma\in \Gamma\}$. 

If $h\in H^2$, $f\in H^\infty$ and $\gamma\in \Gamma$, then 
\[
M_fC_\gamma(h) = f(h\circ \gamma) = (f\circ \gamma^{-1}\gamma )(h\circ \gamma) = ((f\circ \gamma^{-1})h)\circ \gamma = C_\gamma M_{f\circ \gamma^{-1}}(h).
\]
Hence, $M_{f\circ \gamma^{-1}} = C_\gamma^{-1}M_fC_\gamma$. This shows that the action of $\Gamma$ on $H^\infty$ is implemented by a similarity. The fixed point space for this action is the algebra $H^\infty_\Gamma$ which is the set of $f\in H^\infty$ such that $M_fC_\gamma = C_\gamma M_f$ for all $\gamma\in \Gamma$. We now give the most basic connection between $\alg(S,\Gamma)$ and $H^\infty_\Gamma$.

\begin{prop}
The commutant of $\alg(S,\Gamma)$ in $B(H^2)$ is $H^\infty_\Gamma$.
\end{prop}

\begin{proof}
If $T\in \alg(S,\Gamma)'$ , then $TS=ST$ which forces $T = M_f\in H^\infty$. Since $M_fC_\gamma= C_\gamma M_f$, for all $\gamma\in \Gamma$, we get $f\in H^\infty_\Gamma$. 
\end{proof}

\subsection{Inner multipliers of the defect space.}
While the space $N$ does not contain analytic functions, it was shown by Forelli~\cite[Lemma~7]{forelli} that it was possible to multiply $N$ into $H^2$ by an inner function. The following observation about the invariant subspaces of $\alg(S,\Gamma)$ will prove useful.

\begin{prop}\label{alginvariant}
The invariant subspaces of $\alg(S,\Gamma)$ are of the form $\phi H^2$, where $\phi$ is a character-automorphic inner.
\end{prop}

\begin{proof}
If $\cal{M}$ is invariant for $\alg(S,\Gamma)$, then $\cal{M}$ is shift invariant. It follows that $\cal{M} = \phi H^2$ for some inner function $\phi$. Since $\cal{M}$ is invariant under $C_\gamma$ we get that $(\phi \circ \gamma) H^2 = \phi H^2$. It follows, from the uniqueness statement in Beurling's theorem~\cite[Chapter~4.4]{helson}, that $\phi \circ \gamma = \sigma(\gamma) \phi$ where $\sigma(\gamma) \in \bb{T}$.  If $\gamma_1,\gamma_2\in \Gamma$, then 
\[
\sigma(\gamma_1\gamma_2)\phi =\phi\circ (\gamma_1\gamma_2) 
=(\phi \circ \gamma_1)\circ \gamma_2 
 =(\sigma(\gamma_1)\phi)\circ \gamma_2
=\sigma(\gamma_1)\sigma(\gamma_2)\phi.
\]
Hence, $\sigma\in \hat{\Gamma}$. 
\end{proof}

Our next proposition shows that there are no finite-dimensional subalgebras of $H^\infty$, aside from the constants. We will need this fact for the special case where the subalgebra under consideration is $H^\infty_\Gamma$. 

\begin{prop}\label{infdim}
If $\cal{A}$ is a unital subalgebra of $H^\infty$ that contains nonconstant functions, then $\cal{A}$ is infinite dimensional. 
\end{prop}

\begin{proof}
Assume to the contrary that $\cal{A}$ is $n$-dimensional. Let $f\in \cal{A}$ be nonconstant. By subtracting $f(0)$, we may assume that $f(0)=0$ and $f\not = 0$. The elements $1,f,f^2,\ldots,f^n$ must be linearly dependent and so there exists $a_0,\ldots,a_n$ such that $\sum_{j=0}^n a_j f^j =0$. Evaluating at $z=0$, we get $a_0=0$ and so $f(a_1+a_2f+\ldots+a_nf^{n-1})=0$. Since $f\not = 0$, $a_1+a_2f+\ldots+a_nf^{n-1}=0$. Repeating the above argument yields $a_1=\ldots=a_n=0$, a contradiction. 
\end{proof}

As noted by Forelli, the subspace $H^2_\Gamma+N$ is invariant for $H^\infty_\Gamma$. To see this pick $h\in H^2_\Gamma+N$ and choose $k\in H^2$ such that $\Psi(k) = h$. If $f\in \cl{H^2_{0,\Gamma}}$ and $g\in H^\infty_\Gamma$, then 
\[
\int f\cl{gh} = \int f\cl{g}\cl{\Psi(k)} = \int f\cl{gk} = 0.
\]
Hence, $H^\infty_\Gamma(H^2_\Gamma+N)\subseteq H^2_\Gamma+N$. This shows that $N$ is semi-invariant for $H^\infty_\Gamma$ and so the compression of $H^\infty_\Gamma$ to $N$ is a homomorphism of $H^\infty_\Gamma$.  Since $H^\infty_\Gamma$ is infinite dimensional and $N$ is finite dimensional we see that this homomorphism has a nontrivial kernel. Hence, there exists a function $f\in H^\infty_\Gamma$ such that $M_f(N) \perp N$. 

While the functions in $N$ are not analytic, they can be multiplied into $H^2$ by a character automorphic inner function. Let $\cal{N}$ be the set of functions $f\in H^2$ such that $fN \subset H^2$. The subspace $\cal{N}$ is closed and invariant for $\alg(S,\Gamma)$. Hence, the subspace $\cal{N}=\phi H^2$ for some character-automorphic inner function $\phi$, by Proposition~\ref{alginvariant}. This means that $\phi N \subseteq H^\infty$ with  $\phi$ character automorphic. This is the proof of~\cite[Lemma~7]{forelli}

We point out that there is a related way to obtain this last fact. Let $\cal{M}$ denote the smallest shift invariant subspace of $L^2$ that contains $H^2_\Gamma+N$. Clearly $\cal{M} = [H^\infty(H^2_\Gamma+N)]$. By the Helson-Lowdenslager theorem~\cite[Theorem~18]{helson}, the subspace is either of the form $\chi_E L^2$ or of the form $\phi H^2$ where $\phi$ is unimodular. Suppose that $\cal{M}=\chi_E L^2$ and note that $\chi_E = 1$, since $1\in \cal{M}$. If $\cal{M}=L^2$, then $\Psi(\cal{M}) = L^2_\Gamma$. A typical element of $\cal{M}$ can be approximated by sums of elements of the form $fg$ where $f\in H^\infty$ and $g\in H^2_\Gamma+N$. Note that $\Psi(fg) = \Psi(f)g$, which shows 
\begin{align*}
L^2_\Gamma&= \Psi(\cal{M}) \subseteq [\{\Psi(fg)\,:\,f\in H^\infty,g\in H^2_\Gamma+N\}]\\
&=[\{\Psi(f)g\,:\,f\in H^\infty,g\in H^2_\Gamma+N\}] \\
&=[\{fg\,:\,f\in H^\infty_\Gamma+N,g\in H^2_\Gamma+N\}]\\
&\subseteq H^2_\Gamma+N+N.N, 
\end{align*}
which is impossible since $N$ is finite dimensional and $\cl{H^\infty_{0,\Gamma}}$ is not. Hence, $\cal{M} = \psi H^2$ for a unimodular character automorphic function $\psi$. Since $1\in H^2_\Gamma+N\subseteq \psi H^2$, we get $1 = \psi\phi$ for an inner function $\phi$ and so $\cl{\psi} = \phi$ is inner. Hence, $\phi(H^2_\Gamma+N)\subseteq H^2$. This ability to multiply $N$ into $H^\infty$ with a character automorphic inner function is central to our proof of Theorem~\ref{distform} and so we record this fact.

\begin{lemma}[Forelli]\label{multN}
There exists a character automorphic inner function $\phi$ such that $\phi N\subseteq H^\infty$.
\end{lemma}

Using Lemma~\ref{multN} we get the following lemma:

\begin{lemma}\label{multk}
If $w\in \bb{D}$ and $\Gamma$ is an admissible Fuchsian group, then there exists a character automorphic inner function $\phi_w$ such that $\phi_w \cl{k_w^\Gamma}\in H^\infty$. 
\end{lemma}

\begin{proof}
We can assume that $H^\infty_\Gamma$ is nontrivial, otherwise $k_w^\Gamma = 1$. Let $\cal{N} = \{f\in H^2\,:\, f\cl{k_w^\Gamma} \in H^2\}$. The subspace $\cal{N}$ is closed subspace and invariant for $\alg(S,\Gamma)$. Since $H^\infty_\Gamma$ is nontrivial, we can find $f\in H^\infty_\Gamma$ such that $f\not =0$ and $f(w)=0$. Let $g\in H^2_\Gamma$ and note that $\inp{f\cl{k_w^\Gamma}}{\cl{g}}= \inp{fg}{k_w^\Gamma} = f(w)g(w)=0$. Hence, $f\cl{k_w^\Gamma}\in H^2_{0,\Gamma}+N$. By Proposition~\ref{multN}, there exists an inner function $\phi$ such that $\phi f\cl{k_w^\Gamma}\in H^2$. It follows that the space $\cal{N}$ is nontrivial and, by Beurling's theorem, there exists an inner function $\phi_w$ such that $\cal{N} = \phi_w H^2$. It follows from Proposition~\ref{alginvariant} that $\phi_w$ is character automorphic and $\phi_w\cl{k_w^\Gamma}\in H^2$. Since $k_w^\Gamma\in H^\infty$, it follows that $\phi_w \cl{k_w^\Gamma}\in H^\infty$.
\end{proof}

\section{A generalization of Abrahamse's theorem}\label{abrahamse}

We now come to a formula for the distance of an element in $L^\infty_\Gamma$ from the algebra $H^\infty_\Gamma$. Actually we will prove a slightly stronger statement. Let $z_1,\ldots,z_n\in \bb{D}$ and let $\cal{I}$ denote the ideal of functions in $H^\infty_\Gamma$ that vanish at the $n$ points $z_1,\ldots,z_n\in \bb{D}$. Let $\cal{K}_\sigma$ denote the span of the kernel functions $k_{z_1}^\sigma,\dots,k_{z_n}^\sigma$ and let \[\cal{N}_\sigma = H^2_\sigma \ominus \cal{K}_\sigma = \{f\in H^2_\sigma\,:\, f(z_1)=\ldots=f(z_n) = 0\}.\]
By combining Theorem~\ref{cyclic} with~\cite[Lemma~4.6]{raghupathi} we get that $\cal{N}_\sigma = [\cal{I}u]$ for any outer function $u\in H^2_\sigma$.

Our first result in this section is a distance formula. The formula is an analogue of the Nehari's theorem which relates the distance of $f\in L^\infty$ from $H^\infty$ to the norm of the Hankel operator with symbol $f$. 

\begin{thm}\label{distform}
Let $\Gamma$ be an admissible group and let $f\in L^\infty_\Gamma$. The distance of $f$ from $\cal{I}$ is given by 
\begin{equation}
\norm{f+\cal{I}}=\sup_{\sigma\in \hat{\Gamma}}\norm{(I-P_{\cal{N}_\sigma})M_fP_{H^2_\sigma}}. 
\end{equation} 
\end{thm}

Here $I$ denotes the identity in $B(L^2)$ and the orthogonal projections are in $B(L^2)$.

\begin{proof}
Let $k_{z_j}^\Gamma$ be the kernel function at the point $z_j$ for the space $H^2_\Gamma$. By duality
\begin{equation}
\label{dist_gamma}\norm{f+\cal{I}}=\sup\mod{\int fg}
\end{equation}
where $g\in \cal{I}_\perp$ and $\norm{g}_1\leq 1$. Since $\cal{I}\subseteq H^\infty_\Gamma$, we see that $\cal{I}_\perp \supseteq H^1_\Gamma+N$. Since the functions in $\cal{I}$ vanish at the points $z_1,\ldots,z_n\in \bb{D}$, we also have $\cl{\cal{K}_\Gamma}\subseteq \cal{I}_\perp$. An application of Proposition~\ref{duals} shows that $\cal{I}_\perp =H^1_{\Gamma,0}+N+\cl{\cal{K}_\Gamma}$. By Proposition~\ref{multk} there exist character automorphic inner functions $\phi_1,\ldots,\phi_n$ such that $\phi_j \cl{k_{z_j}^\Gamma}\in H^2$. Let $\phi_0$ be the inner function, as in Proposition~\ref{multN}, such that $\phi_0 N\subseteq H^\infty$. Let $\phi = \phi_0\ldots\phi_n$ and note that $\phi$ is character automorphic. We have $\phi  g\in H^1$ and, by Proposition~\ref{riesz2}, there exists a character $\sigma\in \hat{\Gamma}$ and an outer function $u\in H^2_\sigma$ such that $\mod{g}=\mod{u}^2$. Rewriting the expression in \eqref{dist_gamma} we get 
\[
\norm{f+\cal{I}}=\sup\mod{\int fu\overline{v}}
\]
where $u$ is an outer function in $H^2_\sigma$ and $\overline{v}u=g$. Since $g\in \cal{I}_\perp$ we see that $\int g h =0$ for all $h\in \cal{I}$. Therefore, $\inp{hu}{v}=0$ and so $v\in L^2\ominus [\cal{I}u]=L^2\ominus \cal{N}_\sigma$. 
Hence, 
\[
\mod{\int fu\overline{v}} = \mod{\inp{fu}{v}}\leq\norm{(I-P_{\cal{N}_\sigma}) M_f P_{H^2_\sigma}}.
\]

The reverse inequality is straightforward since $H^2_\sigma$ is an invariant subspace for the algebra $H^\infty_\Gamma$.
\end{proof}

Our generalization of Abrahamse's theorem now follows from this distance formula. The theorem and its proof both appear to be new. The theorem applies to a fairly broad class of Riemann surfaces. In particular, the theorem applies to quotients of the disk $\bb{D}$ by an action of a finitely generated Fuchsian group. We also point out that the proof, and the lemmas on which it depends, does not require knowledge of function theory on multiply connected domains. 

\begin{thm}Let $\Gamma$ be an admissible Fuchsian group. Let $z_1,\ldots,z_n\in \bb{D}$ and $w_1,\ldots,w_n\in \bb{C}$. There exists a function $f\in H^\infty_\Gamma$ with $\norm{f}_\infty\leq 1$ such that $f(z_j)=w_j$ if and only if the matrices
\begin{equation}
A_\sigma := [(1-w_i\overline{w_j})K^\sigma(z_i,z_j)]_{i,j=1}^n\geq 0 
\end{equation}
for all $\sigma\in \hat{\Gamma}$. 
\end{thm}

\begin{proof}
 If $f\in H^\infty_\Gamma$, then the subspaces $H^2_\sigma$ and $\cal{N}_\sigma$ are invariant for $H^\infty_\Gamma$. The subspace  $\cal{K}_\sigma = H^2_\Gamma\ominus \cal{N}_\sigma$ is semi-invariant for $H^\infty_\Gamma$. Therefore, $\norm{(I-P_{\cal{N}_\sigma})M_fP_{H^2_\sigma}} = \norm{P_{\cal{K}_\sigma}M_f P_{\cal{K}_\sigma}}$. 

Standard facts about multiplier algebras of reproducing kernel
 Hilbert spaces tell us that $\norm{P_{\cal{K}_\sigma}M_f P_{\cal{K}_\sigma}}\leq 1$ if and only if the matrices $A_\sigma\geq 0$ for all $\sigma \in \hat{\Gamma}$, where $w_j = f(z_j)$; see~\cite[Theorem~6.3]{paulsen}.

Assume first that $f\in H^\infty_\Gamma$, with $\norm{f}_\infty\leq 1$ and $f(z_j) = w_j$. Then $ \norm{P_{\cal{K}_\sigma}M_f P_{\cal{K}_\sigma}}\leq \norm{f}_\infty\leq 1$ for all $\sigma\in \hat{\Gamma}$.

For the converse, assume that $A_\sigma\geq 0$ for all $\sigma \in \hat{\Gamma}$. By applying~\cite[Lemma~5.8]{raghupathi} we see that  there exists a function $g\in H^\infty_\Gamma$ such that $g(z_j) = w_j$. Using the distance formula in Proposition~\ref{distform} we get that $\norm{g+\cal{I}} =\sup_{\sigma\in \hat{\Gamma}}\norm{P_{\cal{K}_\sigma}M_f P_{\cal{K}_\sigma}}\leq 1$. For each $n\in \bb{N}$, there exists $h_n\in \cal{I}$ such that $\norm{g+h_n}_\infty\leq 1+n^{-1}\leq 2$. Since $H^\infty_\Gamma$ is $\wk$ closed, the bounded sequence $f_n:=g+h_n$ has a $\wk$ convergent subnet. Let $f$ be the limit of this subnet. Since a point evaluation at a point in the disk is $\wk$ continuous, it follows that $\norm{f}_\infty\leq 1$ and  $f(z_j) = w_j$.
\end{proof}

Abrahamse's original theorem is now the special case where $\Gamma$ is a group of deck transformations. 

\section{Two Examples}\label{examples}

\subsection{Amenable admissible groups}
In this section we consider two specific examples of admissible Fuchsian groups. Our purpose with these examples is to illustrate the concepts outlined in this paper. Our first example arises as the group of deck transformations for a covering map of an annulus $\bb{A}$. The group in our second example is isomorphic to  $\bb{Z}_2\ast \bb{Z}_2$ and the associated fixed-point algebra is the set of ``even'' functions in $H^\infty(\bb{A})$. Note that both $\bb{Z}$ and $\bb{Z}_2\ast \bb{Z}_2$ are amenable. If $\Gamma$ is an amenable group and $s$ is an invariant mean on $\ell^2(\Gamma)$ we can construct a contractive, unital projection $\Phi_s:H^\infty\to H^\infty_\Gamma$. We will use this fact later on. On the other hand if $\Gamma$ is not amenable, then a result of Barrett~\cite{barrett} shows that no such projection can exist.

The first example can also be found in Abrahamse~\cite[Page~298]{abrahamse3} and Sarason~\cite{sarason2}. Let $\bb{A}$ be the annulus with inner radius $r$ and outer radius $R$. The strip $\cal{S}:=\{z\,:\,0<\Re(z)<\pi\}$ is easily seen to be a covering space for the annulus. The covering map is given by $E(z)=\exp\left(\dfrac{z}{\pi}\log\left(\dfrac{r}{R}\right)\right)$. The disk and strip are conformally equivalent and we see that the corresponding group of deck transformations is an infinite cyclic group, which we denote $\Gamma_1$, generated by $\gamma(z)=\dfrac{z-a}{1-az}$, where $0<a<1$. 

The second example is generated by $\beta(z)=-z$ and $\gamma(z)=\dfrac{a-z}{1-az}$ with $a\in (0,1)$. Denote by $\Gamma_2$ the group generated by $\beta$ and $\gamma$. Every element of $\Gamma_2$ can be identified with a word in $\beta$ and $\gamma$. Since $\beta$ and $\gamma$ have order 2, the elements of $\Gamma_2$ are words in $\beta$ and $\gamma$ with the property that every word is an alternating string of $\beta$'s and $\gamma$'s. Hence, there are precisely two words of each length, distinguished by the ``letter'' they begin with. Set $\alpha=\beta\gamma$. Note that $\alpha(z)=\dfrac{z-a}{1-az}$ and so $\Gamma_2$ contains $\Gamma_1$. A word in $\Gamma_2$ is of the form $\alpha^m$ or $\alpha^m\beta$ where $m\in\bb{Z}$. Since $\beta(0)=0$ we see that $\alpha^m(0)= \alpha^m\beta(0)$. If $\alpha^m(0)=0$, then lemma \ref{hyppower} shows $m=0$ and so $\alpha^m$ is not the identity map for $m\not = 0$. It follows that $\alpha^m\beta$ is also not the identity map for $m\not=0$. We have shown that $\Gamma_2\cong \bb{Z}_2\ast \bb{Z}_2$, the free product of $\bb{Z}_2$ with itself, which is an amenable group.

For these two examples, it can be shown directly that the spaces $H^\infty_\Gamma$ are nontrivial, that the associated Blaschke product converges, and that we can compute the corresponding character $\sigma_0$. Note that the elements of $H^\infty_{\Gamma_2}$ are exactly the functions in $H^\infty_{\Gamma_1}$ that satisfy $f(z)=f(\beta(z))=f(-z)$. We will call such a function \textit{even}. The proof of the following lemma involves an an elementary induction argument.

\begin{lemma}\label{hyppower}
Let $\gamma(z)=\dfrac{z-a}{1-az}$ where $a\in (-1,1)$ and set 
\[
a_n=\dfrac{(1+a)^n-(1-a)^n}{(1-a)^n+(1+a)^n} 
\]
for $n\geq 1$. We have $\gamma^{(n)}(z)=\dfrac{z-a_n}{1-a_nz}$ and $\gamma^{(-n)}(z)=\dfrac{z+a_n}{1+a_nz}$ for $n\geq 1$.
\end{lemma}

From Lemma~\ref{hyppower} we see that $\Gamma_1$ is admissible.

\begin{lemma}
The Blaschke sum $\sum\limits_{\gamma\in\Gamma_1} (1-\mod{\gamma(0)})$ converges. 
\end{lemma}

\begin{proof}
We have that the elements of $\Gamma_1$ are of the form $\gamma^{(n)}$ where $\gamma(z)=\dfrac{z-a}{1-az}$ with $a\in (0,1)$ and $n\in \bb{Z}$. From the previous lemma we have 
\[
1-a_n=\dfrac{2(1-a)^n}{(1+a)^n+(1-a)^n}\leq 2\left(\dfrac{1-a}{1+a}\right)^n. 
\]
The latter series is geometric and so the above sum converges.
\end{proof}

One consequence of this result is that the Blaschke product with zero set $\gamma^{(n)}(0)$ is convergent. Our next proposition gives us $\sigma_0$.  

\begin{prop}
 Let $B$ denote the Blaschke product associated to the pair $(\Gamma_1,0)$. Then $B\circ \gamma = -B$ and $B\circ \beta = -B$. In particular, $B^2\in H^\infty_{\Gamma_2}$ and the algebra $H^\infty_{\Gamma_2}$ is nontrivial.
\end{prop}

\begin{proof}
The Blaschke product for the orbit of 0 under $\Gamma_1$ is the Blaschke product for the set $\{0,a_n,-a_{n}\,:\,n\geq 1\}$. For $n\in \bb{Z}$, let $\phi_n$ denote the simple Blaschke factor at $\gamma^{(n)}(0)$. For $n\geq 1$, the simple Blaschke factor at $a_n=\gamma^{(-n)}(0)$ is given by 
\[
\dfrac{\mod{a_n}}{a_n}\dfrac{a_n-z}{1-a_nz}=-\dfrac{z-a_n}{1-a_nz}=-\gamma^{(n)}(z).
\]
The factor at $-a_n=\gamma^{(n)}(0)$ can be computed similarly and is 
\[
\dfrac{\mod{-a_n}}{-a_n}\dfrac{-a_n-z}{1+a_nz}=\dfrac{z+a_n}{1+a_nz}=\gamma^{(-n)}(z).
\]
The factor at $0$ is $z$.  Our calculation shows that
\[
\phi_n=
\begin{cases}
\gamma^{(-n)} & n\geq 1 \\ 
z & n=0\\ 
-\gamma^{(-n)} & n\leq -1
\end{cases}. 
\]
Let 
\[
B_N:=\prod_{\mod{j}\leq N}\phi_j =(-1)^N\prod_{\mod{j}\leq N}\gamma^{(j)}. 
\]
Composing with $\gamma$ and multiplying by $\gamma^{(-N)}$ we get
\begin{align*}
\gamma^{(-N)}(B_N\circ \gamma)& = \gamma^{(-N)}(-1)^N\prod_{\mod{j}\leq N}\gamma^{(j+1)} \\
&=((-1)^N\prod_{\mod{j}\leq N}\gamma^{(j)})\gamma^{(N+1)}= B_N \gamma^{(N+1)}.
\end{align*}
Hence, $\gamma^{(-N)}(z)B_N(\gamma(z))=\gamma^{(n+1)}(z)B_N(z)$. Since $a_n\to 1$ as $n\to \infty$ we get that $\gamma^{(n)}(z)\to -1$ and $\gamma^{(-n)}(z)\to 1$ as $n\to \infty$. Hence, on taking the limit in $N$ we get that $B(\gamma(z))=-B(z)$. 

For $n\in \bb{Z}$ we have $\gamma^{(n)}(-z)=-\gamma^{(-n)}(z)$. Hence, $B_N(-z)=(-1)^{2N+1}B_N(z)=-B_N(z)$ and on taking the limit we see $B(-z)=-B(z)$. 
\end{proof}

The stabilizer at the origin for the group $\Gamma_2$ consists of the identity and the automorphism $\beta$. Hence, by Proposition~\ref{blaschkedivision}, the Blaschke product associated with the pair $(\Gamma_2,0)$ is $B^2$. It follows from Proposition~\ref{onb} that $\{B^{2n}\,:\,n\geq 0\}$ is an orthonormal basis for $H^2_{\Gamma_2}$. 

\subsection{Interpolation results}
Since $B^2$ is inner, we see that $H^\infty_{\Gamma_2}$ is really the span of the powers of an inner function. The interpolation theory for spaces generated by a single inner function  is quite simple as Theorem~\ref{interp_inner} will show. To prove Theorem~\ref{interp_inner} we require a few preliminary results.

\begin{lemma}
If $\phi$ is an inner function, then $\overline{\phi(\bb{D})}=\overline{\bb{D}}$.
\end{lemma}

\begin{proof}
The operator of multiplication by $\phi$ on $H^2$ is isometric but not unitary. By the Wold decomposition the spectrum of $M_\phi$ is the closed unit disk. Hence, $\overline{\bb{D}}=\sigma(M_\phi)=\overline{\phi(\bb{D})}$.
\end{proof}

\begin{cor}\label{innercomp}
If $\phi$ is an inner function and $f\in H^\infty$, then $\norm{f}_\infty=\norm{f\circ \phi}_\infty$.  
\end{cor}

\begin{proof}
Since $\phi(\bb{D})\subseteq \bb{D}$, $\norm{f\circ \phi}_\infty\leq \norm{f}_\infty$. If $z\in \bb{D}$, then there exists $\phi(z_n)\in \phi(\bb{D})$ such that $\phi(z_n)\to z$ and so $\mod{f(z)}=\lim_{n\to\infty}\mod{f(\phi(z_n))}\leq \norm{f\circ \phi}_\infty$. 
\end{proof}

\begin{prop}\label{interp_inner}
Let $\phi$ be an inner function such that $\phi(0)=0$. Let $H^2_\phi$ be the closed span in $H^2$ of $\{1,\phi,\phi^2,\ldots\}$ and let $H^\infty_\phi$ be the $\wk$ closure of $\{\phi^n\}_{n\geq 0}$ in $H^\infty$. We have the following: 
\begin{enumerate}
\item The set $\{\phi^n\,:\,n\geq 0\}$ is an orthonormal basis for $H^2_\phi$. 
\item The kernel function for $H^2_\phi$ is $K^\phi(z,w)=\dfrac{1}{1-\phi(z)\overline{\phi(w)}}$.
\item The function $f\in H^\infty_\phi$ if and only if $f\in H^2_\phi\cap H^\infty$ if and only if $f=g\circ \phi$ for some $g\in H^\infty$ with $\norm{g}_\infty=\norm{f}_\infty$.
\item The multiplier algebra of $H^2_\phi$ is $H^\infty_\phi$. 
\item The space $H^2_\phi$ is a complete Nevanlinna-Pick space.
\end{enumerate}
\end{prop}

\begin{proof}
\begin{enumerate}
\item Since the function $\phi$ vanishes at the origin we have that $\phi^n$ is orthonormal and since any element of $H^2_\phi$ can be approximated by a finite linear combination of the form $a_0+a_1\phi+\ldots a_n\phi^n$ we see that it is a basis.
\item We have 
\[
K^\phi(z,w)=\sum_{n=0}^\infty\phi(z)^n\overline{\phi(w)^n}=\dfrac{1}{1-\phi(z)\overline{\phi(w)}} 
\]
\item Let us denote by $H^\infty\circ \phi$ the space of functions of the form $\{g\circ \phi\,:\,g\in H^\infty\}$ and by $\cal{A}$ the $\wk$ closed subalgebra of $H^\infty$ spanned by $\{\phi^n\,:\,n\geq 0\}$. The algebra $H^\infty\circ\phi$ is $\wk$ closed and contains $\cal{A}$. If $f\in H^2_\phi\cap H^\infty$, then $f=\sum_{n=0}^\infty a_n \phi^n$ and so $f=(\sum_{n=0}^\infty a_nz^n)\circ\phi$. The sequence $a_n$ is square summable and so $g=\sum_{n=0}^\infty a_n z^n\in H^2$ and $g\circ \phi=f$. However, by Lemma~\ref{innercomp} we get that $\norm{g}_\infty=\norm{f}_\infty$ and so $g\in H^\infty$. This shows that $H^2_\phi\cap H^\infty\subseteq H^\infty\circ \phi$. If $f\in H^\infty\circ \phi$, then $f=g\circ \phi$ and we can choose a net of polynomials $p_t$ such that $p_t\to g$ in the $\wk$ topology. Composition by $\phi$ is $\wk$ continuous and so $p_t\circ \phi \to f$ in the $\wk$ topology and $\cal{A}=H^\infty\circ\phi$. Finally note, since $H^2_\phi$ is closed in $H^2$, that the space $H^2_\phi\cap H^\infty$ is $\wk$ closed in $H^\infty$ and contains $\phi^n$ for all $n$. This proves that $\cal{A}\subseteq H^2_\phi\cap H^\infty$.
\item Since $1\in H^2_\phi$ and multipliers must be bounded we see that $\mult(H^2_\phi)\subseteq H^2_\phi\cap H^\infty = H^\infty_\phi$. On the other hand, it is clear from $H^\infty_\phi=H^\infty\circ \phi$ that $H^\infty_\phi\subseteq \mult(H^2_\phi)$. The equality of norms follows as in Proposition~\ref{multalg}.
\item To see that $K^\phi$ is a complete Nevanlinna-Pick kernel let us consider $n$ points $z_1,\ldots,z_n$ in the disk and $n$ matrices $W_1,\ldots,W_n$ in $M_k$. Suppose that the matrix $[(I-W_iW_j^\ast)K^\phi(z_i,z_j)]_{i,j=1}^n\geq 0$. There exists $F\in M_k(H^\infty)$ such that $\norm{F}_\infty\leq 1$ and $F(\phi(z_j))=W_j$ by the classical Nevanlinna-Pick theorem. The function $\tilde{F}=F\circ\phi$ has the same norm as $F$ and is in $H^\infty_\phi$. Also $\tilde{F}(z_j)=F(\phi(z_j))=W_j$. 

On the other hand, suppose that there is a function $F$ in $M_k(H^\infty_\phi)$ such that $F(z_j)=W_j$. We know that  $F=G\circ \phi$ for some $G\in M_k(H^\infty)$ with $\norm{G}_\infty=\norm{F}_\infty$. Hence, $G(\phi(z_j))=W_j$ and the corresponding Pick matrix is positive.
\end{enumerate}
\end{proof}

Although the proof of Proposition \ref{interp_inner} is elementary we point out that the representation obtained in part (2) implies the result in part (5) by the work of Agler-McCarthy, McCullough, and Quiggin on complete Nevanlinna-Pick kernels; see Agler and McCarthy's book~\cite{aglermccarthy} for a full description. 

Applying Theorem~\ref{interp_inner} to the space $H^2_{\Gamma_2}$ we get the following matrix-valued interpolation result.

\begin{cor}\label{interpeven}
Given $n$ points $z_1,\ldots,z_n\in \bb{D}$ and $n$ matrices $W_1,\ldots,W_n\in M_k$, there exists a function $F\in M_k(H^\infty_{\Gamma_2})$ with $\norm{F}_\infty \leq 1$ and $F(z_j)=W_j$ if and only if 
\[
\left[\dfrac{I-W_iW_j^\ast}{1-B(z_i)^2\overline{B(z_j)^2}}\right]_{i,j=1}^n = \left[(I-W_iW_j^\ast)K^\Gamma(z_i,z_j)\right]_{i,j=1}^n \geq 0,
\]
where $B$ is the Blaschke product for the orbit $\Gamma_2(0)=\Gamma_1(0)$.
\end{cor}

It is interesting to note that the space $H^\infty_{\Gamma_1}\cong H^\infty(\bb{A})$ contains a subalgebra $H^\infty_{\Gamma_2}$, of ``index 2'', that is the multiplier algebra of a complete Nevanlinna-Pick space. 

We now focus on the fact that both $\Gamma_1$ and $\Gamma_2$ are amenable groups. The following result is well known and we state it for completeness. 

\begin{prop}\label{projection1}
Let $X$ be a Banach space. If $\Gamma$ is an amenable group that acts on $X$, then the map $\Phi$ defined by $\Phi(f)(x)= m(f(\gamma(x)))$ is a linear, contractive, idempotent map whose range is $(X^\ast)_\Gamma$ and $\Phi(\gamma^\ast(f))=\Phi(f)$ for all $\gamma\in \Gamma$. 
\end{prop}

The case we are interested in is $X=L^\infty$. In this case $\Phi$ is a unital positive map and is therefore completely contractive.

If $Y\subseteq X^\ast$ is $\Gamma$-invariant, then so is its $\wk$ closure. To see this note that $Y_\perp$ is $\Gamma$-invariant and so is $\overline{Y}^{\wk}=(Y_\perp)^\perp$. It is a standard argument to check that the above projection maps $H^\infty$ to $H^\infty_\Gamma$. 

In some ways this new projection is an improvement on the one in~\cite{forelli}, since it avoids the defect space $N$, preserves analytic structure and is completely contractive. The following result shows, in the case of $H^\infty$, that $\Phi$ is in some sense the natural projection to consider.

\begin{prop}
Let $f\in H^\infty$ and let the power series expansion of $f\circ \gamma$ be given by $\sum_{n=0}^\infty a_n(\gamma)z^n$. For each $n$, $(a_n(\gamma))_{\gamma\in \Gamma}\in \ell^\infty(\Gamma)$. If $b_n=m(a_n(\gamma))$, then 
\[
\Phi(f)(z)=m(f(\gamma(z)))=\sum_{n=0}^\infty b_nz^n. 
\]
\end{prop}

\begin{proof}
If $k_z$ is the Szeg\"o kernel at $z$, then 
\[
\Gamma(f)(z)=\inp{\Phi(f)}{k_z}=m(\inp{f\circ\gamma}{k_z}=m(f(\gamma(z))). 
\]
We have, $\mod{a_n(\gamma)}=\mod{\inp{f\circ \gamma}{\chi_n}}\leq \norm{f}_\infty$ and so for fixed $n$,  $(a_n(\gamma))_{\gamma\in \Gamma}\in\ell^\infty(\Gamma)$. Now, 
\[
\inp{\Phi(f)}{z^n}=m(\inp{f\circ \gamma}{z^n})=m(a_n(\gamma))=b_n. 
\]
and so $\Phi(f)(z)=\sum_{n=0}^\infty b_n z^n$.
\end{proof}

For the special of $H^\infty$, there is a second averaging result that yields the same projection.

\begin{thm}
Let $\cal{A}\subseteq B(\cal{H})$ and suppose that $\cal{A}=\cal{A}'$. If $\Gamma$ is an amenable group that acts on $\cal{A}$, then there exists a linear, contractive, idempotent map $\Phi:\cal{A}\to \cal{A}_\Gamma$. If $\cal{A}$ is selfadjoint, then so is $\Phi$. 
\end{thm}

\begin{proof}
Let $A\in\cal{A}$ and $f_1,f_2\in \cal{H}$. Define a map $s_A:H\times H\to \bb{C}$ by $s_A(f_1,f_2):=m(\inp{\gamma(A)f_1}{f_2})$ and note this is bounded sesquilinear form on $\cal{H}$ and that $\norm{s_A}\leq \norm{A}$. Therefore, there exists a unique operator $\Phi(A)\in B(H)$ such that $\inp{\Phi(A)f_1}{f_2}=m(\inp{\gamma(A)f_1}{f_2})$. Let $A\in \cal{A}$ and note that 
\begin{align*}
\inp{\Phi(A)Tf_1}{f_2}&=m(\inp{\gamma(A)Tf_1}{f_2})=m(\inp{T\gamma(A)f_1}{f_2})=m(\inp{\gamma(A)f_1}{T^\ast f_2})\\
&=\inp{\Phi(A)f_1}{T^\ast f_2} =\inp{T\Phi(A)f_1}{f_2}.
\end{align*}
Thus, $\Phi(A)\in \cal{A}'=\cal{A}$. The remaining properties are easy to verify. 
\end{proof}

It is easy to check, in the case $\cal{H}=H^2$ and $\cal{A}=H^\infty$, that this agrees with the projection in Proposition~\ref{projection1}. 

\begin{prop}
If $\Phi_j$ is the projection from $H^\infty$ onto $H^\infty_{\Gamma_j}$ for $j=1,2$, then $\Phi_j(A(\bb{D}))$ is the set of constant functions.
\end{prop}

\begin{proof}
Since $\gamma^{(k)}(z)\to -1$ as $k\to \infty$ and $1$ as $k\to -\infty$ we have that, 
\[
m(((\gamma^{(k)}(z))^n)_{k\in \bb{Z}})=\begin{cases}0 & \text{if $k$ is odd} \\ 1 & \text{if $k$ is even} \end{cases}\]
Therefore, $\Phi_j$ maps $z^n$ to a constant and so $A(\bb{D})$ is also mapped to the constant functions. 
\end{proof}

\begin{cor}
There are no nonconstant, $\Gamma_j$-invariant functions in $A(\bb{D})$ for $j=1,2$.
\end{cor}

\begin{proof}
If $f\in A(\bb{D})$ and is fixed by $\Phi_j$, then $f=\Phi_j(f)$ is constant.
\end{proof}

\begin{cor}
The projection $\Phi_j$ is not $\wk$ continuous.
\end{cor}

\begin{proof}
This is a consequence of the fact that the trigonometric polynomials are $\wk$ dense in $H^\infty$ and that $H^\infty_{\Gamma_j}$ contains nonconstant functions.
\end{proof}

Since interpolation is an isometric theory the existence of a contractive projection suggests, at least in the case of amenable groups, that we could approach the problem of interpolation in $H^\infty_\Gamma$ through a related problem on the disk. We will show that while this approach has the positive aspect of providing a fairly simple matrix positivity condition, the condition does not seem refined enough to distinguish the interpolation theory of $H^\infty_{\Gamma_1}$ and $H^\infty_{\Gamma_2}$.

Let $K$ denote the Szeg\"o kernel on the disk.  By a {\itshape weak solution} to the interpolation problem we mean a matrix-valued function $f=\left[f_{i,j}\right]\in M_m(H^\infty)$ such that $\norm{f}_\infty\leq 1$, $f(\alpha(z_l))=W_l$ for all $l=1,\ldots,n$ and $\gamma\in \Gamma$. By a \index{Strong solution}\textit{strong solution} we mean an $f=\left[f_{i,j}\right]\in M_m(H^\infty_\Gamma)$ such that $\norm{f}_\infty\leq 1$ and $f(z_l)=W_l$ for $l=1,\ldots,n$. Clearly every strong solution is a weak solution. If we denote by $\cal{W}$ the set of weak solutions, then it is clear, provided that $\cal{W}\not = \emptyset$, that $\cal{W}$ is a convex, $\wk$-compact subset of the unit ball of $M_m(H^\infty)$.

\begin{thm}
Let $\Gamma$ be an amenable, admissible group. If $f\in M_m(H^\infty)$ is weak solution, then $\Phi_m(f)$ is a strong solution.
\end{thm}

\begin{proof}
We have $\norm{\Phi_m(f)}\leq \norm{f}\leq 1$, $\Phi_m(f)\in M_m(H^\infty_\Gamma)$ and  
\[
\Phi_m(f)(z_l)=(\mu(f_{i,j}(\alpha(z_l))))=W_l. 
\]
\end{proof}

If $Q=\left[q_{i,j}\right]_{i,j\in J}$ is an infinite matrix, then we write $Q\geq 0$ to mean that every finite square submatrix is positive. This is also known as \textit{formal positivity}, but we will suppress the word formal. If $Q$ is already finite, or is the matrix of positive operator on a Hilbert space, then the two notions of positivity coincide. The content of the next result is well known and the proof is a standard \wk{} limit argument.

\begin{prop}\label{inf_NP}
Let $\{z_j\}_{j\in J}$ be a set of points in $\bb{D}$ and $\{W_j\}_{j\in J}\subseteq M_m$. There exists $f\in M_m(H^\infty)$ such that $\norm{f}_\infty\leq 1$ and $f(z_j)=W_j$ if and only the matrix $\left[(I_m-W_iW_j^\ast)K(z_i,z_j)\right]_{i,j\in J}\geq 0$
\end{prop}

If we take the set of points in \ref{inf_NP} to be the orbits of $n$ points $z_1,\ldots,z_n$ under the group $\Gamma_j$, then we obtain the following result:

\begin{cor}\label{inf_NP_2}
Let $\Gamma$ be an amenable, admissible group. Let $P(z,w)$ denote the infinite matrix $\left[K(\gamma(z),\eta(w))\right]_{\gamma,\eta\in \Gamma_j}$. A weak (and hence strong) solution $f\in H^\infty_{\Gamma}$ exists if and only if $\left[(I_m-W_iW_j^\ast)P(z_i,z_j)\right]_{i,j=1}^N\geq 0$
\end{cor}

We know that the condition in Corollary~\ref{inf_NP_2} is equivalent to the condition set forth by Abrahamse for the annulus $\bb{A}$. It would be interesting to have a direct proof that the matrix condition in Corollary~\ref{inf_NP_2} is equivalent to Abrahamse condition that the matrices $A_\lambda$, $\lambda\in \bb{T}$, are all positive semi-definite.

Our examples in this section were primarily designed for illustrative purposes. However, we believe that more constructions of this sort could shed further light on the interpolation problem. 

We close with an observation that relates to the work of McCullough and Paulsen~\cite{mcculloughpaulsen}. Their results show that the $C^\ast$-envelope of the $n$-dimensional quotient algebra $H^\infty(\bb{A})/\cal{I}$ is $M_n(C(\bb{T}))$,	 whenever $n\geq 3$. This should be contrasted with the interpolation result in Corollary~\ref{interpeven} which shows that the $C^\ast$-envelope for the algebra $H^\infty_{\Gamma_2}/\cal{I}$ is $M_n$. It was also shown in~\cite{dprs} that the $C^\ast$-envelope for any two-idempotent operator algebra~\cite{paulsen2} is either $\bb{C}^2$ or $M_2$. In light of this, three aspects of the interpolation problem, in apparently increasing order of difficulty, would seem to be worthy of study. First, the matrix-valued interpolation problem. This would seem to be within reach using the methods outlined in this paper. Second, a more careful examination of the number of characters required to guarantee the existence of a solution to the interpolation problem. Third, a computation of the $C^\ast$-envelope of the quotient algebra $H^\infty_\Gamma/\cal{I}$. In view of the depth of the work in~\cite{fedorovvinnikov} and~\cite{mcculloughpaulsen}, and in view of the results presented here that these latter problems are considerably more challenging.

\end{document}